\documentclass{article}
\usepackage{graphicx} 
\usepackage[utf8]{inputenc}
\usepackage{ragged2e}
\usepackage{amsfonts}
\usepackage{mathtools}
\usepackage{amssymb}
\usepackage{amsthm}
\usepackage{authblk}
\usepackage{amsmath}
\usepackage{multicol}
\usepackage{scalerel}

\numberwithin{equation}{section}
\usepackage{enumitem}
\usepackage{afterpage}
\usepackage{authblk}
\usepackage{float}
\usepackage{tikz-cd}
\usepackage{xcolor}
\usepackage{subcaption}
\usepackage{paracol}
\usepackage{tikz}
\usepackage{listings}
\usetikzlibrary{mindmap}

\usepackage{algorithm}
\usepackage{algorithmic} 
\usepackage{color}
\usepackage{hyperref}
\hypersetup{
    colorlinks=true,
    linktoc=all,    
    linkcolor = blue,
}
\hypersetup{
    colorlinks=false,  
    pdfborder={0 0 0},  
    urlcolor=black      
}

\allowdisplaybreaks

\usepackage{geometry}
\theoremstyle{definition}
\newtheorem{definition}{\bf {Definition}}[section]

\newtheorem{example}[definition]{\bf {Example}}
\newtheorem{corollary}[definition]{\bf {Corollary}}
\newtheorem{lemma}[definition]{\bf {Lemma}}
\newtheorem{theorem}[definition]{\bf {Theorem}}
\newtheorem{conjecture}[definition]{\bf {Conjecture}}

\newtheorem{remark}[definition]{\it Remark}

\newtheorem*{theorem*}{Theorem}

\usepackage{cancel}
\title{
On Some Generalisations of Gauss Sequences}
\author[1,a]{Sathyanarayan Narayan}
\author[1,b]{N. Uday Kiran}
\affil[1]{Department of Mathematics and Computer Science\\
Sri Sathya Sai Institute of Higher Learning, Puttaparthi, India.\thanks{Dedicated to Bhagawan Sri Sathya Sai Baba}}
\affil[a]{\href{mailto:sathyanarayannarayan@sssihl.edu.in}\texttt{sathyanarayannarayan@sssihl.edu.in}}
\affil[b]{\href{mailto:nudaykiran@sssihl.edu.in}\texttt{nudaykiran@sssihl.edu.in}}

\date{\today}

\begin{document}
\maketitle
\begin{abstract}
In this paper, we introduce integer sequences satisfying new congruence properties inspired by the Euler and Gauss congruences, which we call `Euler–Gauss' sequences. Noting that every Gauss sequence is an Euler–Gauss sequence, we compare them with certain generalisations of Gauss sequences and provide several counterexamples. Unlike Gauss sequences, Euler--Gauss sequences include sequences based on distinct prime factors, such as the Smallest Prime Factor and  Greatest Prime Factor sequences (suitably defined at $1$). Moreover, we show that the prime-divisor subclass of Gauss sequences, given by $\sum_{p\mid n}p\,g_p$ for an integer sequence $(g_n)$, admits a natural extension to Euler--Gauss sequences of the form $\sum_{p\mid n}p\,f_p(\operatorname{rad}(n))$, where for each prime $p$, $f_p$ is an integer-valued function and $\operatorname{rad}(n)$ denotes the square-free kernel of $n$. Further, we obtain $q$-analogs of the Euler--Gauss sequences, fill gaps in the literature on $q$-Gauss sequences and also conjecture a divisibility criterion for $q$-Euler--Gauss sequences, which we have verified computationally. We also show that not only do our $q$-Euler--Gauss sequences satisfy the Cyclic Sieving Phenomenon (CSP) exhibited by the $q$-Gauss sequences, but we also derive a new CSP condition for the $\operatorname{SPF}$ and $\operatorname{GPF}$ sequences, not hitherto known in the literature. 
    
    
\end{abstract}
\section{Introduction}\label{sec:1}

It is well-known that many interesting integer sequences are inspired by the Euler congruence
\begin{equation}\label{eq:0}
a^{\varphi(n)} \equiv 1 \pmod{n}, \qquad \gcd(a,n)=1,
\end{equation}
where $\varphi(\cdot)$ denotes the totient function.  
When $n$ is a prime power $p^{r}$, we have $\varphi(p^{r}) = p^{r} - p^{r-1}$, and the congruence simplifies to
$
a^{p^{r}} \equiv a^{p^{r-1}} \pmod{p^{r}},
$
which is satisfied for all integers $a$.  
This identity has inspired a class of sequences, aptly termed \emph{Euler sequences} $(a_n)$, that satisfy the congruence
\begin{equation}\label{eq:001}
a_{p^{r}} \equiv a_{p^{r-1}} \pmod{p^{r}}.
\end{equation}
A natural specialisation of~\eqref{eq:001} leads to the \emph{Gauss sequences}, defined by the congruence
\begin{equation}\label{eq:02}
a_{p^{r}m} \equiv a_{p^{r-1}m} \pmod{p^{r}},
\end{equation}
for all primes $p$ and integers $m,r \ge 1$.  
These sequences are known to exhibit a rich structure and have several equivalent forms. Not only geometric sequences but also several linear recurrence sequences such as the Lucas and Perrin sequences are examples of Gauss sequences.
Curiously, the Fibonacci sequence, despite sharing the same linear recurrence relation as the Lucas numbers, fails to satisfy (\ref{eq:02}). These observations have motivated further study of arithmetic sequences satisfying congruence properties.

In this paper, taking inspiration from the Euler congruence (\ref{eq:0}) and the congruence (\ref{eq:02}), we define and investigate a new class of sequences which we term as the \emph{Euler--Gauss sequences}. This proposed class of sequences indeed satisfies the containment: 
$$
\{\textnormal{Gauss Sequences}\}  \subsetneq \{\textnormal{Euler--Gauss Sequences}\} \subsetneq \{\textnormal{Euler Sequences}\}.
$$
We provide several examples in the paper to distinguish between these classes of sequences. For instance, the Smallest Prime Factor ($\operatorname{SPF}$) and Greatest Prime Factor ($\operatorname{GPF}$) sequences, suitably defined at $1$, are examples of Euler--Gauss sequences but do not satisfy most of the generalisations of the Gauss congruence considered in this paper. In that sense, they are strictly Euler--Gauss sequences. 

We also introduce $q$-analogs of the Euler--Gauss sequences, which we term as \emph{$q$-Euler--Gauss sequences}, and explore their connections to $q$-analogs of other congruence-based sequences. In the course of this exploration, we also suggest modifications to certain well-known sequences. To the best of our knowledge, the classes of sequences discussed in this paper do not appear in the literature on integer sequences and, we believe, may help bridge a fundamental gap between Gauss sequences and Euler sequences. 

The congruence \eqref{eq:02} is popularly known by its equivalent form (see \cite[Exercise 5.2(a)]{stanley_2}),  
\begin{equation}\label{eq:03}
    \sum_{d\mid n}\mu(d)a_\frac{n}{d}\equiv0\pmod{n}
\end{equation}
called the \textit{Gauss congruence}. Although observed first by Gauss that (\ref{eq:03}) holds for geometric sequences, other mathematicians working independently also arrived at this result (see for instance, \cite{dickson,Smyth,zarelua}). Gauss sequences appear by different names in the contemporary literature: 
Gauss sequences in Minton's work \cite{minton_gauss}, Dold sequences in \cite{dold_survey,dold}, and Generalised Fermat sequences and Newton sequences in \cite{generalized_fermat, fern_gauss}. Gorodetsky \cite{goro_q-gauss} generalised the Gauss sequences to satisfy \eqref{eq:03} with respect to a set. Further, W\'ojcik \cite{Wojcik_mobius_replace_gauss} showed that the M\"obius function in \eqref{eq:03} may be replaced by any integer-valued arithmetic function $f$ satisfying $f(1)=\pm1$ and $\sum_{d\mid n}f(d)\equiv0\pmod n$ for all $n$; see also the survey \cite{dold_survey}. 

Several interesting examples of Gauss sequences are found in the literature. Certain constant-term sequences, such as the Apéry sequence, as discussed by Straub in \cite{armin_constant}, are Gauss sequences. All Newton sequences generated by integers, such as the sequence of the sum of divisors of a number, are examples of Gauss sequences (see for instance, \cite{generalized_fermat}). Sequences satisfying the Dwork congruences for all primes (see \cite{dwork_congruence}) are Gauss sequences; for instance, the Apéry sequence. Minton \cite{minton_gauss} characterises all linear recurrence Gauss sequences as rational multiples of trace sequences, which also explains why the Lucas sequence is a Gauss sequence while the Fibonacci sequence is not. Other examples of Gauss sequences include the sequence of traces of powers of integer matrices, $(\text{Tr}(M^n))$.
See also the survey article \cite{zarelua} for further examples.

Among the general divisor-sum representation of Gauss sequences, $a_n=\sum_{d\mid n}dg_d$, where $(g_n)$ is any integer sequence, we consider the case in which the coefficients corresponding to composite divisors are zero, giving the special subclass of sequences of the form $a_n=\sum_{p\mid n}p\,g_p,$
where the sum is over the prime divisors $p$ of $n$. This prime-divisor case admits a natural extension to Euler--Gauss sequences of the form
\begin{equation}\label{eq:029}
a_n=\sum_{p\mid n}p\,f_p(\operatorname{rad}(n)),
\end{equation}
where for each prime $p$, $f_p:\mathbb{Z}_{>0}\to\mathbb{Z}$ is an integer-valued function and $\operatorname{rad}(n)$ denotes the square-free kernel of $n$. This class of sequences forms a $\mathbb{Z}$-module under pointwise addition and scalar multiplication and is also closed under Hadamard product. These sequences and their generalisations are discussed in detail
in Section \ref{sec:4}.

\begin{figure}[t]
\centering
\includegraphics[scale=0.4]{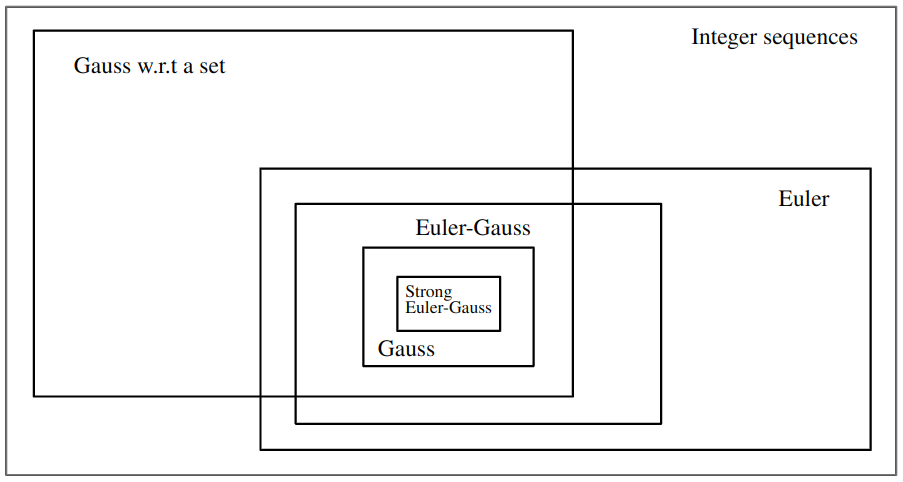}
\caption{Gauss sequences and generalisations}
\label{fig:1}
\end{figure}

We now briefly outline the intuition behind our proposed sequence. 
To this end, let us reconsider the Euler congruence (\ref{eq:0}) under the condition $\gcd(a,n)=1$ and $a\neq 0$. Substituting the Dirichlet convolution representation of the totient function,
\[
\varphi(n)=\sum_{d\mid n}\mu(d)\,\frac{n}{d},
\]
in (\ref{eq:0}) and rewriting the resulting expression in product form, we obtain
\[
a^{\varphi(n)}=a^{\sum_{d|n}\mu(d)n/d}=\prod_{d\mid n} \bigl(a^{\,n/d}\bigr)^{\mu(d)} \equiv 1 \pmod{n}.
\]
Using $\gcd(a,n)=1$ and segregating the factors by $\mu(d)$, this simplifies to the congruence
$$
\prod_{\substack{d\mid n \\ \mu(d)=1}} a^{\frac{n}{d}}\equiv \prod_{\substack{d\mid n \\ \mu(d)=-1}} a^{\frac{n}{d}}\pmod{n}.$$ 
Motivated by the above congruence, we now define our Euler--Gauss congruence for integer sequences $(a_n)$ as follows:
\begin{equation}\label{eq:04}
    A_n^+\equiv A_n^-\pmod{n},
\end{equation}
where 
\[A_n^+=\prod_{\substack{d\mid n \\ \mu(d)=1}} a_{\frac{n}{d}}\text{,\hspace{4mm}and\hspace{4mm}}A_n^-=\prod_{\substack{d\mid n \\ \mu(d)=-1}} a_{\frac{n}{d}}\text{.}\]
We call the sequences $(a_n)$ satisfying (\ref{eq:04}) as the Euler--Gauss sequences. It is noteworthy that the Gauss sequences satisfying (\ref{eq:03}) can also be expressed in a similar fashion as (\ref{eq:04}) albeit in summation form. Although Gauss sequences are defined with a summation, they are indeed Euler--Gauss sequences (see Theorem \ref{thm:4}). The Venn diagram in Figure \ref{fig:1} illustrates containment relations between several classes of sequences. The terminology in the diagram is explained in Section \ref{sec:1.1} and Section \ref{sec:2}.

Based on the Gauss congruence (\ref{eq:03}), Gorodetsky \cite{goro_q-gauss} recently introduced a $q$-analog of Gauss sequences known as the $q$-Gauss sequences. These sequences $(a_{n}(q))$ with terms in $\mathbb{Z}[q]$ satisfy the congruence
\begin{equation}
  \sum_{d\mid n}\mu(d)a_\frac{n}{d}(q^d)\equiv0\pmod{[n]_q}, 
\end{equation}
for all $n\geq1$, where $[n]_q=1+q+\cdots+q^{n-1}$ denotes the $q$-analog of the integer $n$. Further, he discussed the $q$-Gauss sequences in connection with the Cyclic Sieving Phenomenon (CSP). 
He showed that for a $q$-Gauss sequence $(a_n(q))$ with non-negative integer coefficients and a sequence of sets $(X_n)$ with $|X_n|=a_n(1)$ and the addition modulo $n$ group $\mathbb{Z}_n$ acting on $X_n$ in such a way that 
\begin{equation}\label{eq:020}
    |X_{n}^j|=|X_{\gcd(n,j)}|
\end{equation}
for all $n\geq1$, where $X_n^j$ denotes the fixed point set of $j\in\mathbb{Z}_n$ acting on $X_n$, the triple $(X_n, \mathbb{Z}_n, a_n(q))$ satisfies the CSP for every $n\geq1$. For a combinatorial interpretation of this CSP, refer to the recent work of Gossow \cite{fern_gauss}. 
  
Motivated by the works of Gorodetsky, we introduce the $q$-Euler--Gauss sequences and explore their connections to other $q$-analogs. We show that our class of sequences naturally fits into various $q$-analogs (see Figure \ref{fig:2}). While establishing connections and presenting counterexamples for different classes in the $q$-analogs, we also identified and addressed some gaps in the current literature. One such gap is addressed in Section \ref{sec:3}, where we introduce a \textit{Modified $q$-Gauss congruence} that admits a characterisation. In contrast, the 
$q$-analog of the congruence (\ref{eq:02}) defined by Gorodetsky (see \cite[Remark 1.2]{goro_q-gauss}) does not have such a result. 

A salient feature of our Euler--Gauss sequences and its $q$-analog is that the class of sequences defined by \ref{eq:029} contains non-trivial examples such as the $k$\textsuperscript{th} Smallest Prime Factor ($k$-$\operatorname{SPF}$) and the $k$\textsuperscript{th} Greatest Prime Factor ($k$-$\operatorname{GPF}$) sequences. These sequences have been studied by Erd\"os, Alladi and several other prominent mathematicians (see \cite{gpfspf0,gpfspf1} and the references therein). The duality identities \cite[Lemma 1]{gpfspf0} proved by Alladi appear for these sequences in the context of Euler--Gauss sequences. Moreover, the CSP for these sequences is new in the literature. For instance, we prove the following result in Section \ref{sec:5} (Corollary \ref{thm:3}) for the $\operatorname{SPF}$ sequence (similar result holds for the $\operatorname{GPF}$ sequence): 

\begin{theorem*}[CSP for the $q$-$\operatorname{SPF}$ sequence]
    Let $(S_n(q))$ be the sequence of $q$-analogs of the $\operatorname{SPF}$ sequence terms. Define $X_{1}=\emptyset$. For an $n>1$, let $X_{n}$ be a set with $|X_{\operatorname{SPF}(n)}|=\operatorname{SPF}(n)$ and $\mathbb{Z}_n$ act on $X_n$ such that  
    \begin{equation}\label{eq:019}        |X_n^j|=\left|X_{\gcd{(\operatorname{SPF}(n), j)}}\right|
    \end{equation} 
     for all $j\in\mathbb{Z}_n$, where $\operatorname{SPF}(n)$ is the smallest prime factor of $n$. 
     Then $\left(X_n,\mathbb{Z}_n,S_n(q)\right)$ is a CSP triple.  
\end{theorem*}

By comparing \eqref{eq:019} and \eqref{eq:020}, this CSP condition can be easily contrasted with the standard CSP for $q$-Gauss sequences. We further show in Section \ref{sec:5} (Corollary \ref{coro:1}) that the standard \emph{CSP} extends to the $q$-Euler--Gauss sequences $(a_{n}(q))$ as well, for the triples $(X_{n}, \mathbb{Z}_{n}, a_{n}(q))$, where \eqref{eq:020} holds when $a_d(1)\neq0$ for all proper divisors $d\mid n$ such that $\frac{n}{d}$ is not a prime power. These observations suggest that $q$-Euler--Gauss sequences may admit a combinatorial interpretation.  

The rest of the paper is organised as follows. In the next section, some important preliminaries are discussed, and a concise background on Gauss sequences is provided. Following this, in Section \ref{sec:2}, the Euler--Gauss sequences are explored through examples and containment results are proved. In Section \ref{sec:3}, the $q$-Euler--Gauss sequences are introduced and compared with generalisations of the $q$-Gauss sequences in the literature. Further, a necessary condition for sequences in $\mathbb{Z}[q]$ to be $q$-Euler--Gauss sequences is determined and stated in the form of a divisibility criterion verified computationally up to $n=10^8$ terms (See Conjecture \ref{conj:1}). Subsequently, in Section \ref{sec:4}, these results are applied, and the class of sequences characterised by \eqref{eq:029} is discussed. In particular, the $\operatorname{SPF}$ and $\operatorname{GPF}$ sequences are explored in detail. In Section \ref{sec:5}, the CSP for $q$-Gauss sequences is extended to the $q$-Euler--Gauss sequences and the new CSP condition for the $\operatorname{SPF}$ and $\operatorname{GPF}$ sequences is discussed. 
 
\section{Preliminaries and Background}\label{sec:1.1}
 
In this section, we will discuss some standard generalisations of Gauss sequences from the literature. These sequences are characterised by various modifications of the Gauss congruence. 
We also recall certain equivalents of the Gauss congruence that shall be used to establish containment relations for the Euler--Gauss sequences. To this end, some terminology is introduced for certain equivalents that shall appear frequently in the following sections.

In this paper, we restrict our attention only to integer sequences. 
It is important to note that Minton \cite{minton_gauss} defines the Euler and Gauss sequences for rational numbers and prefixes the integer cases with \emph{strict}. We will not follow this practice. The works by Gorodetsky \cite{goro_q-gauss}, Gossow \cite{fern_gauss}, and Minton \cite{minton_gauss} are key references for this paper and contain important concepts, generalisations and classification results that constitute a background for this work.

\begin{definition}[Euler Sequence, \cite{minton_gauss}]
    An integer sequence $(a_n)$ is defined to be an Euler sequence if 
    \begin{equation}\label{eq:06}
        a_{p^r}\equiv a_{p^{r-1}}\pmod{p^r},
    \end{equation}
    for all primes $p$ and integers $r\geq1$. 
\end{definition}

As discussed in the introduction section, (\ref{eq:03}) defines the Gauss sequences. Although the conditions (\ref{eq:02}) and (\ref{eq:03}) are equivalent \cite[Exercise 5.2(a)]{stanley_2}, we employ separate definitions for clarity of exposition. 
\begin{definition}[Gauss Congruence and Pairwise Gauss Congruence]
We say that an integer sequence $(a_n)$ satisfies the Pairwise Gauss congruence with respect to a prime $p$ if  
    \begin{equation}\label{eq:07}
        a_{p^rm}\equiv a_{p^{r-1}m}\pmod{p^r},
    \end{equation}
for all integers $r\geq1$ and $m\geq1$. The sequence $(a_{n})$ is said to satisfy the Gauss congruence if 
\begin{equation}\label{eq:071}
    \sum_{d\mid n}\mu(d)a_\frac{n}{d}\equiv0\pmod{n}\text{ .}
\end{equation}
\end{definition}

When \eqref{eq:07} is satisfied for all primes, it is equivalent to the Gauss congruence. Since \eqref{eq:06} is a special case of \eqref{eq:07}, this implies that every Gauss sequence is an Euler sequence. But the converse is not true. 
 
For example, the sequence
\[a_n=\begin{cases}
    1&n\text{ is a prime power}\\
    0&\text{otherwise},
\end{cases}\] is an Euler sequence but not a Gauss sequence.

Gorodetsky \cite{goro_q-gauss} generalises the Gauss sequences to satisfy the Gauss congruence with respect to a set. The set here is a non-empty subset of the set of all primes and will be denoted by $S$. The set of all primes will be denoted by $\mathbb{P}$.

\begin{definition}[Gauss sequence with respect to $S$, \cite{goro_q-gauss}]\label{def:1} 
    An integer sequence $(a_n)$ is defined to be a Gauss sequence with respect to $S$ if  
    \begin{equation}\label{eq:072}
      \sum_{d\mid n}\mu(d)a_\frac{mn}{d}\equiv0\pmod{n}  
    \end{equation}
    for all integers $n\geq1$ divisible only by primes from $S$ and all integers $m\geq1$.  
\end{definition}

Definition \ref{def:1} appears to be inspired by the Pairwise Gauss congruence (\ref{eq:07}). Indeed, \cite[Lemma 2.1]{goro_q-gauss} shows that a sequence $(a_n)$ is a Gauss sequence with respect to $S$ if and only if it satisfies the Pairwise Gauss congruence for all primes in $S$. Clearly, if $S=\mathbb{P}$, (\ref{eq:072}) is equivalent to (\ref{eq:071}). In contrast, if $S=\{2\}$ then the sequence $(2n)$ is a Gauss sequence with respect to $S$ but not a Gauss sequence.

The characterisation of Gauss sequences $(a_n)$ by the divisor sum 
\begin{equation}\label{eq:09}
    a_n=\sum_{d\mid n}dg_d
\end{equation} arises by applying M\"{o}bius inversion to the Gauss congruence. A generalisation of this divisor sum forms a subclass of the Euler--Gauss sequences, which we shall discuss in a later section. To the best of our knowledge, it has not been noted in the literature that Gauss sequences can also be characterised via Lambert series. 
A Lambert series is defined by a sequence $(b_{n})$ as a series of the form
$
\sum_{n\geq 1} b_n q^{n}/(1-q^{n}).
$
For a formal power series $\sum_{n\geq 1} a_n q^n$ arising from a sequence $(a_n)$, we have
\[
\sum_{n\geq 1} a_n q^n = \sum_{n\geq 1} b_n \frac{q^n}{1 - q^n},
\]
and the coefficient sequences $(a_n)$ and $(b_n)$ are related by
\begin{equation}\label{eq:021}
    a_n = \sum_{d \mid n} b_d .
\end{equation}
Comparing the terms of \eqref{eq:021} and \eqref{eq:09}, one can see that if $(a_n)$ is a Gauss sequence, then the sums in \eqref{eq:09} and \eqref{eq:021} become equal. Further, by a M\"{o}bius inversion, since both $(g_n)$ and $(b_n)$ are uniquely determined by $(a_n)$, $b_{n}=ng_{n}$ for all $n\geq 1$. We therefore obtain the following characterising condition for Gauss sequences, expressed in terms of Lambert series.  

\begin{theorem*}
 An integer sequence $(a_n)$ is a Gauss sequence if and only if its generating function has a Lambert series expansion with integer coefficients $b_n$ such that $n\mid b_n$ for all $n\geq 1$. 
\end{theorem*}

Stanley's book \cite{stanley_2} contains a characterising condition for the Gauss sequences in terms of their generating functions. An integer sequence $(a_n)$ is a Gauss sequence if and only if
\begin{equation}\label{eq:08}
\exp\left({\sum_{n\geq1}\frac{a_n}{n}q^n}\right)\in\mathbb{Z}[[q]]\text{ , }
\end{equation}
where $\mathbb{Z}[[q]]$ denotes the ring of formal integer power series. Gorodetsky extends this characterisation to the Gauss sequences with respect to the singleton set $\{p\}$, where $p$ is prime. He showed that $(a_n)$ is a Gauss sequence with respect to $\{p\}$ if and only if 
\[\exp\left(\sum_{n\geq1}\frac{a_n}{n}q^n\right)\in\mathbb{Z}_p[[q]]\text{ , }\] where $\mathbb{Z}_p[[q]]$ denotes the ring of formal power series with coefficients as $p$-adic integers. 
 
In \cite{minton_gauss}, Minton characterises the linear recurrence Gauss sequences. He determines that an integer linear recurrence sequence $(a_n)$ is a rational multiple of a Gauss sequence if and only if it is a \textit{trace} sequence. A sequence of rational numbers $(a_n)$ such that, 
    \[a_n=\sum_{j=1}^{r}\alpha_j\text{Tr}_{\frac{K}{\mathbb{Q}}}(\theta_j^n)\text{ , }\]
where $K$ is an algebraic number field, $\alpha_j \in \mathbb{Q}$ and $\theta_j \in K$ for all $1\leq j\leq r$, is a trace sequence. For instance, the Lucas sequence $(L_n)$ satisfies $L_n=\phi_+^n+\phi_-^n=\text{Tr}_\frac{\mathbb{Q}\sqrt{5}}{\mathbb{Q}}(\phi_+^n)$ and is a trace sequence. Here, $\phi_+$ and $\phi_-$ are the roots of the characteristic polynomial $q^2-q-1$ associated with the linear recurrence relation for the Lucas sequence, $L_n=L_{n-1}+L_{n-2}\text{, }n\geq2$. Although sharing the same linear recurrence relation, the Fibonacci sequence $(F_n)$ satisfies $F_n=\phi_+^n-\phi_-^n$ and is not a trace sequence. 

Polynomial recurrence Gauss sequences are also discussed in the literature. For instance, Straub \cite{armin_constant} gives an important sequence that originated as a `constant term', namely the sequence of Apéry numbers $(a_n)$ defined by the polynomial recurrence relation, 
 \[(n+1)^3a_{n+1}-(34n^3+51n^2+27n+5)a_n+n^3a_{n-1} = 0, \\\text{\hspace{5mm}}a_0=1, a_1=5.\]
 Moreover, he proved that the class of constant term sequences $(a_n)$ defined by 
\[a_n=\text{ct}[P(\mathbf{x})^n]\text{ , }\]
where $P(\mathbf{x})$ is a multivariate Laurent polynomial in $\mathbf{x}=(x_1,\cdots,x_d)$ and $\text{ct}[\cdot]$ denotes `constant term of', are Gauss sequences (see \cite[Proposition 5.1]{armin_constant}).

\subsection{\texorpdfstring{$q$}{q}-Analogs}
Gorodetsky \cite{goro_q-gauss} and Gossow \cite{fern_gauss} discuss interesting equivalents and generalisations of the $q$-Gauss sequences. In this subsection, some key congruences that characterise these generalisations are discussed. Some equivalents that we arrived at by using simpler techniques are explored. Gaps in the definitions of $q$-analogs of certain congruences are identified, and some terminology is introduced to compare them with our modified $q$-analogs. An introduction to the CSP for the $q$-Gauss sequences is also provided.
\begin{definition}[$q$-Gauss sequence, \cite{goro_q-gauss}]
    A sequence $(a_n(q))$ in $\mathbb{Z}[q]$ is defined to be a $q$-Gauss sequence if  it satisfies the congruence 
    \begin{equation}\label{eq:05}
        \sum_{d\mid n}\mu(d)a_\frac{n}{d}(q^d)\equiv0\pmod{[n]_q}
    \end{equation}
    for all $n\geq1$. The congruence \eqref{eq:05} is known as the $q$-Gauss congruence.  
    
\end{definition}
Gorodetsky extended the Pairwise Gauss congruence to a $q$-analog setting. Although he defined it as in \eqref{eq:010}, he did not attach a name to the congruence.  For clarity of exposition, we use the following terminology.
\begin{definition}[Weak $q$-Gauss Sequence, \cite{goro_q-gauss}]
    For a sequence $(a_n(q))$ in $\mathbb{Z}[q]$ and a fixed prime $p$, the congruence 
    \begin{equation}\label{eq:010}
        a_{p^rm}(q)\equiv a_{p^{r-1}m}(q^p)\pmod{[p^r]_q}
    \end{equation}
    where integers $m,r\geq1$, is defined as the $q$-analog of the Pairwise Gauss congruence with respect to the prime $p$. A sequence $(a_n(q))$ that satisfies \eqref{eq:010} for all primes $p$ is defined to be a weak $q$-Gauss sequence.     
\end{definition}

The $q$-analog differs from the integer case in the sense that the $q$-Gauss congruence is not equivalent to the $q$-analog of the Pairwise Gauss congruence when satisfied for all primes. Indeed, every $q$-Gauss sequence satisfies \eqref{eq:010} at all primes $p$ (see \cite[Lemma 2.2]{goro_q-gauss}), but the converse need not hold. See Section \ref{sec:3} for further details. In this section, we introduce a modification of the $q$-analog of the Pairwise Gauss congruence so that we have an equivalence to the $q$-Gauss congruence when satisfied for all primes $p$.   

Gorodetsky extended the Gauss congruence with respect to a set to the $q$-analog setting. The set is a non-empty subset $S\subseteq\mathbb{P}$.
\begin{definition}[$q$-Gauss Sequence with respect to $S$, \cite{goro_q-gauss}]
    A sequence $(a_n(q))$ in $\mathbb{Z}[q]$ is a $q$-Gauss sequence with respect to $S$ if
    \begin{equation}\label{eq:017}
        \sum_{d\mid n}\mu(d)a_\frac{nm}{d}(q^d)\equiv 0\pmod{[n]_q}
    \end{equation}
    for all integers $n\geq1$ divisible only by primes from $S$ and all integers $m\geq1$.  
\end{definition}

It was established in \cite{goro_q-gauss} that a $q$-Gauss sequence with respect to $S$ satisfies the $q$-analog of the Pairwise Gauss congruence for all primes in $S$. However, the converse is not true. Therefore, we again suggest a modification of the $q$-Gauss congruence with respect to $S$ and introduce a new congruence that is equivalent to our modified $q$-analog of the Pairwise Gauss congruence when satisfied for all primes in $S$ (see Section \ref{sec:3}).

As a special case of Gorodetsky's result \cite[Proposition 2.5]{goro_q-gauss} it is true that for every $q$-Gauss sequence $(a_n(q))$, there exists an integer sequence $(g_n)$ such that,   
\[
    a_n(q)\equiv\sum_{d\mid n, d<n}\left[\frac{n}{d}\right]_{q^d}g_\frac{n}{d}\pmod{[n]_q}
\]
for all $n$.  

Nevertheless, as in the integer setting, simply applying M\"{o}bius inversion to the definition of the $q$-Gauss congruence leads us to the following result.
\begin{theorem*}
    A sequence $(a_n(q))$ in $\mathbb{Z}[q]$ is a $q$-Gauss sequence if and only if there exists a sequence $(g_n(q))$ in $\mathbb{Z}[q]$ such that
\begin{equation}\label{eq:022}
a_n(q)=\sum_{d\mid n}\left[\frac{n}{d}\right]_{q^d}g_{\frac{n}{d}}(q^d)
\end{equation}
for every $n\geq1$.
\end{theorem*}
\begin{proof}
For the sequence $(a_n(q))$, the definition of the $q$-Gauss congruence implies that there exists a sequence $(g_n(q))$ in $\mathbb{Z}[q]$ such that
\[
\sum_{d\mid n}\mu(d)a_{\frac{n}{d}}(q^d)=g_n(q)[n]_q
\]
for every $n\geq1$. Substituting this identity into \eqref{eq:022} and reindexing the resulting double sum by setting $e=dd'$ yields
\[
\sum_{d\mid n}\left[\frac{n}{d}\right]_{q^d}g_{\frac{n}{d}}(q^d)
=\sum_{d\mid n}\sum_{d'\mid\frac{n}{d}}
a_{\frac{n}{dd'}}(q^{dd'})\mu(d')
=\sum_{e\mid n}a_{\frac{n}{e}}(q^e)\sum_{d'\mid e}\mu(d')
=a_n(q),
\]
where the last equality follows from the identity
$\sum_{d\mid m}\mu(d)=\delta_{m,1}$. Hence, \eqref{eq:022} holds.

Conversely, suppose that $a_n(q)$ is defined by \eqref{eq:022}. Then, again setting $e=dd'$,
\[
\sum_{d\mid n}\mu(d)a_{\frac{n}{d}}(q^d)
=\sum_{d\mid n}\mu(d)\sum_{d'\mid\frac{n}{d}}
\left[\frac{n}{dd'}\right]_{q^{dd'}}
g_{\frac{n}{dd'}}(q^{dd'})
=\sum_{e\mid n}\left[\frac{n}{e}\right]_{q^e}
g_{\frac{n}{e}}(q^e)\sum_{d'\mid e}\mu(d')
=g_n(q)[n]_q.
\]
Therefore, $[n]_q$ divides
$\sum_{d\mid n}\mu(d)a_{\frac{n}{d}}(q^d)$ and the $q$-Gauss congruence holds for $n$. Since $n$ was arbitrary, $(a_n(q))$ is a $q$-Gauss sequence.
\end{proof}
Although the above characterisation result is not recorded in the current literature, we use it extensively in Section \ref{sec:3} to arrive at conditions for $q$-Euler--Gauss sequences to satisfy the $q$-Gauss congruence.  
Since the polynomials $\left[\frac{n}{d}\right]_{q^d}$ in \eqref{eq:022} satisfy
\begin{equation}\label{eq:023}
    \left[\frac{n}{d}\right]_{\omega_n^{jd}}=
    \begin{cases}
        \frac{n}{d}&\frac{n}{d}\mid j\\
        0&\text{otherwise}
    \end{cases}
\end{equation}
for all integers $n,j\geq1$ and $n$\textsuperscript{th} primitive roots of unity $\omega_n$, this result is useful for computation. In fact, the following key characterising condition for $q$-Gauss sequences is a direct consequence of this result. Since the expansion \eqref{eq:022} for $q$-Gauss sequences did not feature in Gorodetsky's arguments to arrive at this condition, for the sake of completeness, we supply a proof of the same with \eqref{eq:022}. 

\begin{theorem*}[\protect{\cite[Corollary 2.3]{goro_q-gauss}}]
A sequence $(a_n(q))$ is a $q$-Gauss sequence if and only if it satisfies
\begin{equation}\label{eq:011}
    a_n(\omega_n^j)=a_{\gcd(n,j)}(1)
\end{equation}
for all integers $n,j\geq1$ and $n$\textsuperscript{th} primitive roots of unity $\omega_n$.
\end{theorem*}
\begin{proof}
    For a $q$-Gauss sequence $(a_n(q))$, the characterisation \eqref{eq:022} holds so that for each divisor $d\mid n$, 
    \[a_n(\omega_n^j)=\sum_{\substack{d\mid n\\\frac{n}{d}\mid j}}\frac{n}{d}g_\frac{n}{d}(1)=\sum_{d\mid\gcd(n,j)}dg_d(1)=a_{\gcd(n,j)}(1)\] 
    and \eqref{eq:011} is satisfied. The first equality is due to the property \eqref{eq:023} of the polynomials $[\frac{n}{d}]_{q^d}$. For the converse, since it is true that for any prime $p$, $a_p(\omega_p^j)-a_1(1)=0$ for all $i$ such that $p\nmid j$, it follows that $a_p(q)=[p]_qg_p(q)+g_1(q^p)$ for some $g_p(q),g_1(q)\in\mathbb{Z}[q]$ with $a_1(q)=g_1(q)$. With this as the base case for induction, we assume that for all divisors $d\mid n$, $d<n$, the characterisation \eqref{eq:022} holds. Then,  
    \[f(q)=a_n(q)-\sum_{d\mid n, d>1}\left[\frac{n}{d}\right]_{q^d}g_\frac{n}{d}(q^d)\]
    and at roots of unity $\omega_n^j$ where $n\nmid j$, it evaluates to
    \[f(\omega_n^j)=a_{\gcd(n,j)}(1)-\sum_{\substack{d\mid n\\\frac{n}{d}\mid j}}\frac{n}{d}g_\frac{n}{d}(1)=a_{\gcd(n,j)}(1)-\sum_{d\mid\gcd(n,j)}dg_d(1)=0\text{ .}\]
    Therefore, it must be divisible by $[n]_q$ and the characterisation \eqref{eq:022} holds for $a_n(q)$. 
\end{proof}
The characterising condition (\ref{eq:011}) shows that $q$-Gauss sequences possess enumerative behaviour. Using this condition, Gorodetsky \cite{goro_q-gauss} and Gossow (\cite{fern_gauss}) establish a link between $q$-Gauss sequences and the Cyclic Sieving Phenomenon(CSP) as discussed in Reiner, Stanton and White \cite{csp}. Gorodetsky \cite{goro_q-gauss} showed that for every $q$-Gauss sequence, a sequence of triples exhibiting the CSP can be generated. Each triple consists of $a_n(q)$ from a $q$-Gauss sequence with non-negative coefficients $(a_n(q))$, a finite set $X_n$ with cardinality $|X_n|=a_n(1)$, and $\mathbb{Z}_n$ acting on $X_n$ such that
\[|X_n^j| = |X_{\gcd(n,j)}|\text{ .}\] 
 Each such triple $(X_n, \mathbb{Z}_n, a_n(q))$ satisfies the CSP since, 
\[|X_n^j|=|X_{\gcd(n,j)}|=a_{\gcd(n,j)}(1)=a_n(\omega_n^j)\]
where the last equality is due to \eqref{eq:011}.

Some interesting examples of the CSP obtained using $q$-Gauss sequences are available in Gossow's work \cite{fern_gauss}. He observes that every Gauss sequence is simply an enumerator for a sequence of sets satisfying the fixed point counting condition above, calling such a sequence of sets a `Lyndon structure'. This is inspired by the work \cite{csp-lyndon} on `Lyndon-like cyclic sieving' by Alexandersson, Linusson, and Potka.
 
\section{Euler--Gauss Sequences}\label{sec:2} 
In this section, the class of Euler--Gauss sequences is compared with generalisations of the Gauss sequences in the literature. The class of `strong Euler--Gauss sequences', inspired by Euler's congruence \eqref{eq:0}, is introduced. The containment relationships represented in Figure \ref{fig:1} are established, and a condition for the Euler--Gauss sequences to become Gauss sequences is determined. Further, strict examples of these classes are discussed in detail with the required theoretical results to support them.  

\begin{definition}[Euler--Gauss Sequence]
    An integer sequence $(a_n)$ is defined to be an Euler--Gauss sequence if it satisfies 
    \begin{equation}\label{eq:012}
        A_n^+\equiv A_n^-\pmod{n}
    \end{equation}
    for all $n\geq1$, where 
    \[A_n^+=\prod_{\substack{d\mid n \\ \mu(d)=1}} a_{\frac{n}{d}}\text{,\hspace{4mm}and\hspace{4mm}}A_n^-=\prod_{\substack{d\mid n \\ \mu(d)=-1}} a_{\frac{n}{d}}\text{ .}\] We term the congruence \eqref{eq:012} as the `Euler--Gauss congruence'.
\end{definition}

From this definition, one can see that the Euler--Gauss sequences are closed under the Hadamard product. This property can be used to generate several examples of these sequences. One can also see that the Euler--Gauss sequences are indeed Euler sequences.
Indeed, at prime powers $p^r$, since the M\"{o}bius function is non-zero only for the divisors $1$ and $p$, the products on either side of the Euler--Gauss congruence simplify to
\[A_{p^r}^+=a_{p^r}\text{, \hspace{4mm}and \hspace{4mm}}A_{p^r}^-=a_{p^{r-1}}\text{ .}\]
Then, by definition the Euler--Gauss congruence becomes $a_{p^r}\equiv a_{p^{r-1}}\pmod{p^r}$ which precisely defines the Euler sequences. 
 
\begin{theorem}\label{thm:4}
    Every Gauss sequence is an Euler--Gauss sequence.
\end{theorem}
\begin{proof}
    Let $(a_n)$ be a Gauss sequence. Then, for $n=p^rm$, $p^{r+1}\nmid n$, it follows that, 
    \begin{align*}
        \prod_{\substack{d\mid n\\\mu(d)=1}}a_\frac{n}{d}
        &=\prod_{\substack{d\mid m\\\mu(d)=1}}a_\frac{p^rm}{d}\cdot\prod_{\substack{d\mid m\\\mu(d)=-1}}a_\frac{p^{r-1}m}{d}\\
        &\equiv \prod_{\substack{d\mid m\\\mu(d)=1}}a_\frac{p^{r-1}m}{d}\cdot\prod_{\substack{d\mid m\\\mu(d)=-1}}a_\frac{p^{r}m}{d} \pmod{p^r}\\
        &\equiv \prod_{\substack{d\mid n\\\mu(d)=-1}}a_\frac{n}{d} \pmod{p^r} \text{ .}
    \end{align*}
Since this is true for all prime divisors of $n$, the Euler--Gauss congruence is satisfied.
\end{proof}

It is noteworthy that if we replace \eqref{eq:012} with the stronger requirement $(A_n^+)(A_n^-)^{-1}\equiv1\pmod{n},$
then any sequence satisfying the congruence for all $n\geq1$ must have all its terms in $\{1,-1\}$. Indeed, fix $n\geq1$ and set $m=np$ where $p$ is any prime. Since $\mu(p)=-1$ and $a_n=a_{m/p}$, $a_n$ occurs as a factor of $A_m^-$. Now, for the inverse of $A_m^-$ modulo $m$ to exist, we must have $\gcd(A_m^-,m)=1$, and hence $\gcd(a_n,p)=1$. As $p$ was an arbitrary prime, $a_n$ must indeed be $1$ or $-1$. 

At this point, one can either simply relax the coprimality condition or, instead, compute the inverse modulo $n$ only after cancelling the common factors in $(A_n^+)(A_n^-)^{-1}$. In the latter case, the product becomes
\[
\left(\frac{A_n^+}{\gcd(A_n^+,A_n^-)}\right)
\left(\frac{A_n^-}{\gcd(A_n^+,A_n^-)}\right)^{-1}.
\]
and the congruence takes the form
\begin{equation}\label{eq:013}
    \left(\frac{A_n^-}{\gcd(A_n^+,A_n^-)}\right)^{-1}\equiv\left(\frac{A_n^+}{\gcd(A_n^+,A_n^-)}\right)^{-1}\pmod{n}\text{ .}
\end{equation}
This congruence can be satisfied non-trivially. For instance, if $A_n^-\mid A_n^+$ for all $n$, then, $\prod_{d\mid n}a_\frac{n}{d}^{\mu(d)}$ is an integer for all $n$. The congruence \eqref{eq:013} becomes equivalent to, 
\[\prod_{d\mid n}a_\frac{n}{d}^{\mu(d)}=nb_n+1\]
for some integer sequence $(b_n)$. Further, by applying M\"{o}bius inversion for the product case, we obtain 
\begin{equation}\label{eq:014}
    a_n=\prod_{d\mid n}(db_d+1)\text{.}
\end{equation}
For all integer sequences $(a_n)$ of this form, the congruence \eqref{eq:013} is indeed non-trivially satisfied. Thus, the sequences characterised by this congruence deserve investigation. By definition, as the Euler--Gauss congruence necessarily holds for these sequences, we term them as `strong Euler--Gauss sequences'.
\begin{definition}[Strong Euler--Gauss Sequence]\label{defn:strongEG}
    An integer sequence $(a_n)$ is defined to be a strong Euler--Gauss sequence if for all $n\geq1$  
    \begin{enumerate}[label=\roman*.]
        \item $\frac{A_n^+}{\gcd(A_n^+,A_n^-)}$ and $\frac{A_n^-}{\gcd(A_n^+,A_n^-)}$ are each coprime to $n$ so that their inverses modulo $n$ exist and, 
        \item $(a_n)$ satisfies the congruence \eqref{eq:013}
    \end{enumerate}
We term the congruence \eqref{eq:013} as the strong Euler--Gauss congruence. 
\end{definition}
\begin{remark}\label{rem:01}
    Assuming the coprimality condition (i.), the strong Euler--Gauss congruence \eqref{eq:013} is equivalently expressed as
    \[\frac{A_n^+}{A_n^-}\equiv 1\pmod{n}\]
    where the fraction $A_n^+/A_n^-$ is reduced to lowest terms, i.e., after cancellation of the common factors. 
\end{remark}

Note that, in all computations in this paper, we adopt the convention $0^0 = 1$. Thus, in the above definitions, for every factor $a_{\frac{n}{d}}$ with $\mu(d)=0$, we have $a_{\frac{n}{d}}^{\mu(d)} = 1$. 

We now show that every strong Euler--Gauss sequence is a Gauss sequence. To this end, a pairwise congruence for the strong Euler--Gauss sequences is obtained. 
\begin{lemma} \label{lem:01}
    An integer sequence $(a_n)$ is a strong Euler--Gauss sequence if and only if
    \begin{enumerate}
        \item[i.] $\frac{a_{p^rm}}{\gcd\left(a_{p^rm},a_{p^{r-1}m}\right)}$ and $\frac{a_{p^{r-1}m}}{\gcd\left(a_{p^rm},a_{p^{r-1}m}\right)}$ are each coprime to $p$ and, 
        \item[ii.] $(a_n)$ satisfies the congruence 
        \begin{equation}\label{eq:015}
        \left(\frac{a_{p^rm}}{\gcd(a_{p^rm},a_{p^{r-1}m})}\right)^{-1}\equiv\left(\frac{a_{p^{r-1}m}}{\gcd(a_{p^rm},a_{p^{r-1}m})}\right)^{-1}\pmod{p^r}
        \end{equation}
    \end{enumerate}
    for all primes $p$ and integers $r\geq1$ and $m\geq1$.
\end{lemma}
\begin{proof}
    For a strong Euler--Gauss sequence $(a_n)$, it follows by definition that, for $n=p^r$, the coprimality condition (i.) and the congruence \eqref{eq:015} imply the Pairwise Gauss congruence \eqref{eq:07}, and in particular the congruence \eqref{eq:06}, where $p$ is prime and $r\geq1$. Equivalently, every strong Euler--Gauss sequence is an Euler sequence. Thus, for the forward implication, an induction argument with this as the base case may be used.
    
    For $n=p^rm$, $r,m\geq1$, write 
    $A_{p^rm}^+=a_{p^rm}P^+$ and $A_{p^rm}^-=a_{p^{r-1}m}P^-$
    where \[P^+=\prod_{\substack{d\mid m, d>1\\\mu(d)=1}} a_\frac{p^rm}{d}\prod_{\substack{d\mid m \\\mu(d)=-1}} a_\frac{p^{r-1}m}{d}\hspace{4mm}\text{ and }\hspace{4mm}P^-=\prod_{\substack{d\mid m \\\mu(d)=-1}} a_\frac{p^{r}m}{d}\prod_{\substack{d\mid m, d>1\\\mu(d)=1}} a_\frac{p^{r-1}m}{d}.\]
    Assuming the coprimality condition (i.) for every divisor $d \mid n$, $d<n$, the integers $a_{\frac{p^rm}{d}}/\gcd(a_{\frac{p^rm}{d}}, a_{\frac{p^{r-1}m}{d}})$ and  $a_{\frac{p^{r-1}m}{d}}/\gcd(a_{\frac{p^rm}{d}}, a_{\frac{p^{r-1}m}{d}})$ are then each coprime to $p$ so that the exponents of $p$ in the prime factorisations of $a_{p^rm}$ and $a_{p^{r-1}m}$ are equal. Consequently, the same holds for $P^+$ and $P^-$ and further since it holds for $A_{p^rm}^+$ and $A_{p^rm}^-$ by the definition of a strong Euler--Gauss sequence,  the exponents of $p$ in the prime factorisation $a_{p^rm}$ and $a_{p^{r-1}m}$ are equal. Equivalently, $a_{p^rm}/\gcd(a_{p^rm}, a_{p^{r-1}m})$ and  $a_{p^{r-1}m}/\gcd(a_{p^rm}, a_{p^{r-1}m})$ are each coprime to $p$. 

    If the congruence \eqref{eq:015} is assumed to hold for all $d\mid n$, $d<n$, then \[\frac{a_{p^rm/d}}{a_{p^{r-1}m/d}}\equiv1\pmod{p^r}\] where the fraction is reduced to lowest terms. Hence, $P^+/P^-\equiv1\pmod{p^r}$, where the fraction is again reduced to lowest terms. By Remark \ref{rem:01},  $A_{p^rm}^+/A_{p^rm}^-\equiv1\pmod{p^r}$ so that the factors $P^+$, $P^-$ cancel, yielding $\frac{a_{p^rm}}{a_{p^{r-1}m}}\equiv1\pmod{p^r}$ where the fraction is reduced to lowest terms.  

    For the reverse implication, assume that the coprimality condition (i.) and the congruence \eqref{eq:015} hold. Let $n$ be written as $p^rm$ for a prime $p$ and integers $r,m\geq1$. Then, by (i.), the exponents of $p$ in the prime factorisations of  $A_{p^rm}^+=a_{p^rm}P^+$ and $A_{p^rm}^-=a_{p^{r-1}m}P^-$ are equal, and by (ii.), \[\frac{a_{p^rm}P^+}{a_{p^{r-1}m}P^-}\equiv1\pmod{p^r}\] where the fraction is reduced to lowest terms. Since this holds for all primes $p$ occuring in the prime factorisation of $n$, the conditions for $(a_n)$ to be a strong Euler--Gauss sequence are satisfied.
\end{proof} 
\begin{theorem}\label{thm:strongEG}
    Every strong Euler–Gauss sequence is a Gauss sequence. 
\end{theorem}

\begin{proof}
    From Lemma \ref{lem:01}, a strong Euler--Gauss sequence $(a_n)$
    satisfies \[\frac{a_{p^rm}}{\gcd(a_{p^rm},a_{p^{r-1}m})}\equiv\frac{a_{p^{r-1}m}}{\gcd(a_{p^rm},a_{p^{r-1}m})}\pmod{p^r}\]
    for all primes $p$ and $r\geq1$,$m\geq1$. This is true since \eqref{eq:015} is the same as this congruence, but with the added assumption of coprimality. Further, the $\gcd(a_{p^rm}, a_{p^{r-1}m})$ factor may be multiplied to either side of this congruence, and it will still hold. The resulting congruence $a_{p^rm}\equiv a_{p^{r-1}m}\pmod{p^r}$ is precisely the Pairwise Gauss congruence \eqref{eq:07} with respect to $p$. Since this is true for all $p$, the desired containment is proved.
 
\end{proof}
We now show that by simply assuming that the terms of an Euler--Gauss sequence $(a_n)$ are coprime to $n$, they satisfy the Gauss congruence property. An analogue of this result is established for the $q$-Euler--Gauss sequences in the next section (see Theorem~\ref{thm:1}). 
\begin{theorem}
    Let $(a_n)$ be an Euler--Gauss sequence such that $\gcd(a_n,n)=1$ for all $n\geq1$. Then $(a_n)$ is a Gauss sequence.
\end{theorem}
\begin{proof}
    We prove the result by induction. For every prime power $p^r$, the Euler--Gauss congruence coincides with the Pairwise Gauss congruence, and hence the base case holds.

    Writing $n=p^r m$ where $p$ is a prime divisor of $n$, so that $r\geq1$ and $\gcd(p,m)=1$, we assume that the Pairwise Gauss congruence holds for every proper divisor of $n$. Then, the Euler--Gauss congruence at $n$ yields
    \[ a_{p^rm}\prod_{\substack{d\mid m,\,d>1\\\mu(d)=1}}
        a_\frac{p^rm}{d}
        \prod_{\substack{d\mid m\\\mu(d)=-1}}
        a_\frac{p^{r-1}m}{d}
        \equiv a_{p^{r-1}m}\prod_{\substack{d\mid m,\,d>1\\\mu(d)=1}}
        a_\frac{p^{r-1}m}{d}
        \prod_{\substack{d\mid m\\\mu(d)=-1}}
        a_\frac{p^rm}{d}\pmod{p^r}.
    \]
    By the induction hypothesis, $a_{\frac{p^{r-1}m}{d}}\equiv a_{\frac{p^rm}{d}}\pmod{p^r}$ for every $d\mid m$ with $d>1$ so that
    \begin{align*}
        p^r&\mid
        (a_{p^rm}-a_{p^{r-1}m})
        \prod_{\substack{d\mid m,\,d>1\\\mu(d)=1}}
        a_\frac{p^rm}{d}
        \prod_{\substack{d\mid m\\\mu(d)=-1}}
        a_\frac{p^rm}{d}\text{ or, }\\
        p^r&\mid
        (a_{p^rm}-a_{p^{r-1}m})
        \prod_{\substack{d\mid m,\,d>1\\\mu(d)\neq0}}
        a_\frac{p^rm}{d}.
    \end{align*}
    By $\gcd(a_k,k)=1$, each factor in the product is coprime to $p$ so that the Pairwise Gauss congruence holds at $n=p^rm$, completing the induction. 
\end{proof}
\subsection{Examples} 

In this subsection, we provide some examples of the Euler--Gauss sequences. We also provide counterexamples to show that the containments in Figure \ref{fig:1} are strict. From Theorem \ref{thm:strongEG}, we know that every strong Euler--Gauss sequence is a Gauss sequence. However, the converse does not hold. For example, the Lucas sequence, although a Gauss sequence, fails to satisfy the strong Euler--Gauss congruence for infinitely many values of $n$.

\begin{example}
    The sequence $(a_n)$ defined by 
\[a_n=
    2^{k+1}-1\text{ for } 2^k\mid n, 2^{k+1}\nmid n, k\geq0\]
is a strong Euler--Gauss sequence. That it satisfies the strong Euler--Gauss congruence follows from
\[A_{2^km}^+=\prod_{\substack{d\mid m \\ \mu(d)=1}} (2^{k+1}-1)\cdot\prod_{\substack{d\mid m \\ \mu(d)=-1}} (2^{k}-1)=\prod_{\substack{d\mid m \\ \mu(d)=-1}} (2^{k+1}-1)\cdot\prod_{\substack{d\mid m \\ \mu(d)=1}} (2^{k}-1)=A_{2^km}^-\text{ ,}\]
where $m>1$. This is true since the number of divisors $d$ for any integer $m$, with $\mu(d)=-1$, is the same as the number of divisors $d'$ with $\mu(d')=1$. For powers of $2$, $A_{2^k}^+=2^{k+1}-1$ and $A_{2^k}^-=2^{k}-1$ and therefore their inverses modulo $2^k$ are also congruent.

\end{example}

In the above example the product $\prod_{d\mid n}a_\frac{n}{d}^{\mu(d)}$ is not an integer for some $n$. In case we wish to choose a sequence such that $\prod_{d\mid n}a_\frac{n}{d}^{\mu(d)}$ is an integer, from (\ref{eq:014}) we can design a sequence 
$$
    a_n=\prod_{d\mid n}(db_d+1)\text{,}
$$
for an arbitrary sequence $(b_{n})$. It is also easy to see that this strong Euler--Gauss sequence is indeed a Gauss sequence from the following observation
    \begin{align*}
    a_{p^r m}
    &= \prod_{d\mid m}(db_d+1) \cdot \prod_{d\mid m}(pdb_{pd}+1) \cdots \prod_{d\mid m}(p^r db_{p^r d}+1) \\
    &\equiv \prod_{d\mid m}(db_d+1) \cdots \prod_{d\mid m}(p^{r-1} db_{p^{r-1} d}+1) \pmod{p^r} \\
    &\equiv a_{p^{r-1} m} \pmod{p^r}
\end{align*}
for all primes $p$ and integers $r\geq1$ and $m\geq1$.

\begin{example}\label{ex:2}
    Consider the sequence $(a_n)$ defined as, 
\[
    a_n=
    \begin{cases}
        4 & 12 \nmid n\\
        6(2^{k-2}-1)+1 & 2^k3\mid n\text{, }2^{k+1}\nmid n\text{ }\forall k\geq2
    \end{cases}
    \text{ .}
\]
It is an Euler--Gauss sequence but not a Gauss sequence. This can be verified case-by-case. Clearly, for all $n$ such that $12\nmid n$, the subsequence is constant; thus, it is Euler--Gauss. We consider the other case $n=2^k3^jm$ where $k\geq2$, $j\geq1$ and $m>1$ if $j=1$. For such $n$, by definition,
\begin{align*}
    A_n^+&=\prod_{\substack{d\mid m\\\mu(d)=1}}a_\frac{2^k3^jm}{d}a_\frac{2^{k-1}3^{j-1}m}{d}\prod_{\substack{d\mid m\\\mu(d)=-1}}a_\frac{2^{k-1}3^jm}{d}a_\frac{2^k3^{j-1}m}{d}\\&=\prod_{\substack{d\mid m\\\mu(d)=-1}}a_\frac{2^k3^jm}{d}a_\frac{2^{k-1}3^{j-1}m}{d}\prod_{\substack{d\mid m\\\mu(d)=1}}a_\frac{2^{k-1}3^jm}{d}a_\frac{2^k3^{j-1}m}{d}=A_n^-
\end{align*}
and for the remaining $n=2^k3$ where $k\geq2$,  
\[A_{2^k3}^+-A_{2^k3}^-=\begin{cases}
    12 &k=2\\
    2^{k}3&k>2
\end{cases}
\text{ .}
\]
Therefore, the Euler--Gauss congruence is satisfied in each case. However,  
\[a_{2^k3}-a_{2^{k-1}3}=\begin{cases}
    -3 &k=2\\
    2^{k-2}3&k>2
\end{cases}\]
which is not divisible by $2^k$ in any case, and so it does not satisfy the Gauss congruence with respect to $\{2\}$. Therefore, this sequence is not a Gauss sequence, and it does not satisfy the Pairwise Gauss congruence with respect to $2$.  

For primes $p>2$, $a_{p^rm}=a_{p^{r-1}m}$ for all $r\geq1$ and $m\geq1$ by definition, since it is true that  $2^k3\mid p^rm$ implies $2^k3\mid p^{r-1}m$ if $p\neq2$. Therefore, this sequence is a Gauss sequence with respect to odd primes, $S=\{3,5,7\cdots p \cdots\}$. Further, from the equivalence condition with the Lambert series, the generating function of $(a_n)$ is a Lambert series with the coefficients $(b_n)$ defined as,  
\[
    b_n=
    \begin{cases}
        4 & n=1\\
        -3 & n=12\\
        2^{k-2}3 & n=2^k3, k\geq3\\
        0 & \text{otherwise}
    \end{cases}
\]
This is true as $b_n=\sum_{d\mid n}\mu(d)a_\frac{n}{d}$ and since for all primes $p>2$, $a_{p^rm}=a_{p^{r-1}m}$ for all $r\geq1$ and $m\geq1$, $b_n=0$ for all $n>1$ other than $n=2^k3, k\geq2$. For $n=2^k3$, since 
$b_n=a_{2^k3}+a_{2^{k-1}}-a_{2^{k-1}3}-a_{2^k}$ and $a_{2^{k-1}}=a_{2^{k}}$, the values for $b_n=a_{2^k3}-a_{2^{k-1}3}$ are already computed above and they match. At all $n=2^k3, k\geq2$, clearly $n\nmid b_n$ and, as expected, the Gauss congruence is not satisfied. For other values of $n$, $b_n=0$ and the Gauss congruence is satisfied.

From the equivalence condition with generating functions, the power series $\exp\left(\sum_{n\geq1}\frac{a_n}{n}q^n\right)$ of the sequence does not belong to $\mathbb{Z}[[q]]$. Nevertheless, one can determine that it belongs to $\tilde{S}^{-1}\mathbb{Z}[[q]]$, the formal power series ring with coefficients from the localisation $\tilde{S}^{-1}\mathbb{Z}$, where $\tilde{S}$ denotes the minimal multiplicative set generated by $S=\{2\}$, that is $\tilde{S}=\{1,2,4,8,\cdots,2^k,\cdots\}$.
This can be seen using the relationship 
\begin{equation}\label{eq:027}
\exp\left(\sum_{m\geq1}\frac{a_m}{m}q^m\right)=\prod_{n\geq1}\frac{1}{(1-q^n)^{c_n}},
\end{equation}
where $c_n=\frac{1}{n}\sum_{d\mid n}\mu(d)a_\frac{n}{d}$ for all $n\geq1$. By substitution, $c_n=\frac{b_n}{n}$ where the $b_n$ are the Lambert series coefficients of the generating function of $(a_n)$. Therefore, 
\[c_n=\begin{cases}
    4&n=1\\
    -\frac{1}{4}&n=12\\
    \frac{1}{4}&n=2^k3, k\geq3\\
    0&\text{otherwise}
\end{cases}\text{ .}\]
From this, the expansion of $\exp\left(\sum_{n\geq1}\frac{a_n}{n}q^n\right)$ must be a power series with coefficients whose denominators are powers of $2$, thus belonging to $\tilde{S}^{-1}\mathbb{Z}[[q]]$.
\end{example}
\begin{example}
Another sequence similar to the above example is: 
\[
    a_n=
    \begin{cases}
        6 & 6\nmid n, 15\nmid n\\
        2^k-1 & 2^k3\mid n, 2^{k+1}\nmid n, k\geq 1\\
        1+15\cdot\frac{3^{k-1}-1}{2} & 3^k5\mid n, 3^{k+1}\nmid n, 6\nmid n, k\geq1.
    \end{cases}
\]
Again, this is an Euler--Gauss sequence that does not satisfy the Gauss congruence with respect to the set $\{2,3\}$. For this sequence, the power series $\exp(\sum_{n\geq1}\frac{a_n}{n}q^n)$ belongs to $\tilde{S}^{-1}\mathbb{Z}[[q]]$, the formal power series ring with coefficients from $\tilde{S}^{-1}\mathbb{Z}$, the localisation of the integers at the minimal multiplicative set generated by $S=\{2,3\}$.     
\end{example}
The following example illustrates that a part of the intersection of Euler sequences and Gauss sequences with respect to some $S$ lies outside the Euler--Gauss sequences.
\begin{example}\label{ex:4}
    Consider the sequence $(a_n)$ defined as,  

\[a_n=\begin{cases}
    1&\gcd(n,6)=1\text{, or }n=3^k, k\geq0\\
    3&\text{otherwise}
\end{cases}\]
This sequence is an Euler sequence and a Gauss sequence with respect to $S=\{2\}$. However, it is not an Euler--Gauss sequence.
\end{example}

In the above examples, we find that there are Euler--Gauss sequences which are not Gauss sequences but are Gauss sequences with respect to $S$. This leads to the question: \emph{``Does there exist an Euler--Gauss sequence which is not Gauss with respect to any $S$?''} An observation that can be pursued to provide an answer is that the identity sequence $(n)$ satisfies the Euler--Gauss congruence at all $n\geq1$ that are not prime powers but does not satisfy the Gauss congruence at any $n$. To that end, we prove this observation and, following that, provide a result which demonstrates that the sequence $(n)$ may be used to generate several examples of Euler--Gauss sequences that are not Gauss sequences with respect to any $S$. 

\begin{lemma}\label{lem:03}
    The sequence $(n)$ satisfies the Euler--Gauss congruence for all $n$ that are not prime powers. Moreover, $A_n^+ =A_n^-$ for all such $n$. 
\end{lemma}
\begin{proof}
    For the sequence $(n)$, the natural logarithms of the products $A_n^+$ and $A_n^-$ are
    \[\ln{A_n^+}=\sum_{\substack{d\mid n\\\mu(d)=1}}\ln{\frac{n}{d}}\text{\hspace{5mm} and\hspace{5mm}}\ln{A_n^-}=\sum_{\substack{d\mid n\\\mu(d)=-1}}\ln{\frac{n}{d}}\text{ .}\]
    We know that the sum $\sum_{d\mid n}\mu(d)\ln{\frac{n}{d}}$ which is equal to $(\ln{A_n^+}-\ln{A_n^-})$ vanishes at all $n$ that are not prime powers. This is true since the sum simplifies to the Mangoldt function $\Lambda(n)$, which is defined as
    \[\Lambda(n)=\begin{cases}
        \ln{p}&n=p^k\text{ for some prime }p\\0&\text{otherwise} 
    \end{cases}\]
    and is known to satisfy $\sum_{d\nmid n}\Lambda(d)=\ln{n}$ (see \cite[Chapter 2]{antapostol}). Applying M\"{o}bius inversion to this identity precisely proves that the sum equals $0$ at all $n$ that are not prime powers. 
    Therefore at all such $n$, $\ln{A_n^+}=\ln{A_n^-}$ and $A_n^+=A_n^-$ as desired. 
\end{proof}
\begin{theorem}\label{thm:6}
    Let $(a_n)$ be an integer sequence such that $(na_n)$ is an Euler sequence. Then, $(na_n)$ is an Euler--Gauss sequence. Moreover, $(nb_n)$ is a Gauss sequence if and only if $b_{n} =0$ for all $n$. 
\end{theorem}
\begin{proof}
    The products $A_n^+$ and $A_n^-$ for the sequence $(na_n)$ are
    \[A_n^+=\prod_{\substack{d\mid n\\\mu(d)=1}}\frac{n}{d}\prod_{\substack{d\mid n \\ \mu(d)=1}} a_{\frac{n}{d}}\text{ \hspace{5mm}and\hspace{5mm}}A_n^-=\prod_{\substack{d\mid n\\\mu(d)=-1}}\frac{n}{d}\prod_{\substack{d\mid n \\ \mu(d)=-1}} a_\frac{n}{d}\text{ .}\]
    From Lemma \ref{lem:03}, we know that at all $n$ that are not prime powers, the product $\prod_{\substack{d\mid n\\\mu(d)=1}}\frac{n}{d}$ is a common factor to both $A_n^+$ and $A_n^-$. Since this product is a multiple of $n$, $n$ also divides $A_n^+-A_n^-$ and $(na_n)$ satisfies the Euler--Gauss congruence at $n$. For $n$ that are prime powers, the Euler and Euler--Gauss congruences are the same, so that $(na_n)$ is an Euler sequence implies that it is an Euler--Gauss sequence.

    To see that $(nb_n)$ is a Gauss sequence if and only if $b_n=0$ for all $n\geq1$, consider the Pairwise Gauss congruence \eqref{eq:07} which requires $pnb_{pn}\equiv nb_n\pmod{p}$ prime $p$ and $n\geq1$ such that $\gcd(p,n)=1$. This means that $p\mid b_n$ for any arbitrary prime $p$, but this is only possible if $b_n=0$.
\end{proof}

From this theorem, one can observe that $(na_n)$ may be a Gauss sequence with respect to the set $\{p\}$ if and only if $p\mid a_n$ for all $n\geq1$ such that $\gcd(p,n)=1$. This restriction can easily be avoided. The following examples are immediate from the above theorem and illustrate this property.

\begin{example}
    Consider the sequence $(a_n)$ defined as,  
    \[a_n=\begin{cases}
        0&n=1\\
        n^k&n>1
    \end{cases}\]
    for some $k>1$. One may observe that for every $k>1$, this sequence is an Euler sequence. Therefore, from Theorem \ref{thm:6}, $(a_n)$ is an Euler--Gauss sequence since $n\mid a_n$ for all $n\geq1$. Further, for every prime $p$ and all $n$ coprime to $p$, clearly $p\nmid n$ so that $(a_n)$ is not Gauss with respect to any $\{p\}$ and therefore, any $S$.
\end{example}
\begin{example}\label{ex:3}
    The sequence of the sum of primes dividing $n$, with repetitions (see \cite{alladiintegerlog}), also called the integer logarithm of $n$ (\cite[\textit{OEIS A001414}]{oeis}), or the potency of $n$ (see \cite{potencyofn}), is not an example of an Euler--Gauss sequence. In fact, it is not even an Euler sequence. However, the Hadamard product of this sequence with the identity sequence $(n)$:  
\[a_n=n\sum_{p^r\mid\mid n}rp,\quad \textnormal{ where } p \textnormal{ is prime}, \]
is an Euler--Gauss sequence. Further, it is not a Gauss sequence with respect to any set $S$. This is true since for any prime $p$ and integer $r\geq1$, $a_{p^r}=rp^{r+1}$, $a_{p^{r-1}}=(r-1)p^r$ and the Euler congruence is satisfied by $(a_n)$ at arbitrary $p^r$. Therefore, from Theorem \ref{thm:6}, it is an Euler--Gauss sequence. Let $f(n)=\sum_{p^r\mid\mid n}rp$. Then for any prime $p$, the sequence $(a_n)$ fails to satisfy the Pairwise Gauss congruence at all $pn$ such that $\gcd(n,p)=1$ since $a_{pn}=pn(p+f(n))$, $a_n=nf(n)$ so that $a_{pn}$ is divisible by $p$ but $f(n)$ need not be. Therefore, it is not a Gauss sequence with respect to any set $S$.

The Lambert series coefficients associated with $(a_n)$ are 
    \[b_n=\sum_{p^r\|n}\varphi\left(\frac{n}{p^r}\right)( a_{p^{r}}-a_{p^{r-1}}) \]
    for all $n\geq1$, where $p$ is prime and $\varphi(n)$ denotes the totient function. Let  $g(n)$ be an arithmetic function defined such that $g(n)=p$ when $n$ is a prime power $p^r$, $r\geq1$ and $0$ otherwise. Then, $f(n)=(g*\mathbf{1})_n$ where $\mathbf{1}$ denotes the constant function $\mathbf{1}(n)=1$ for all $n$. Then, with $N(\cdot)$ denoting the arithmetic function $N(n)=n$,  
\begin{align*}
    b_n=(a*\mu)_n&=((N\cdot(g*\mathbf{1}))*\mu)_n=((N\cdot g*N)*\mu)_n=(N\cdot g*\varphi)_n
\end{align*}
since the totient function $\varphi=\mu*N$. Further, 
\[(N\cdot g*\varphi)_n=\sum_{p^r\|n}\sum_{i=1}^rp^{i+1}\varphi\left(\frac{n}{p^i}\right)=\sum_{p^r\|n}\varphi\left(\frac{n}{p^r}\right)((r-1)p\varphi(p^r)+p^{r+1})\]
where $(r-1)p\varphi(p^r)+p^{r+1}$ simplifies to $a_{p^{r}}-a_{p^{r-1}}$ as required. 

From this, using the relationship \eqref{eq:027} between the sequence $(a_n)$ and the sequence of coefficients $(c_n)$ defined as $c_n=\frac{b_n}{n}$, it is clear there are infinitely many $n\geq1$ such that $c_n$ is not an integer. This is consistent with the fact that the sequence $(a_n)$ is not Gauss with respect to any set $S$ and thus the associated power series $\exp(\sum_{n\geq1}a_nq^n/n)$
belongs to $\mathbb{Q}[[q]]$.  

This example can be generalised to all sequences of the form 
\[a_n=n\sum_{p\mid n}f_p(n)p\]
where the sum is over prime divisors $p$ of $n$ and $f_m$ are functions $f_m:\mathbb{Z}_{>0}\mapsto\mathbb{Z}$ for all $n\geq1$. That the $(a_n)$'s are Euler--Gauss sequences follows from Theorem \ref{thm:6}, which requires only that they be Euler sequences. Since for any integer $r\geq1$ and prime $p$, $a_{p^r}=p^{r+1}f_p(p^r)$ and $a_{p^{r-1}}=p^rf_p(p^{r-1})$, so that $p^r$ divides both of them, the Euler congruence is satisfied. Example \ref{ex:3} discussed above is clearly a special case of this class of sequences where $f_m(\cdot)\equiv0$ if $m$ is composite and for $m=p$ where $p$ is prime, $f_p(n)=r$ where $r$ is the largest exponent such that $p^r\mid n$.
\end{example}
A natural generalisation of this class of sequences is obtained by removing the factor $n$. In Section \ref{sec:4}, we determine the conditions under which such sequences, i.e., sequences of the form $\sum_{p\mid n}f_p(n)p$, are Euler--Gauss. We also extend this construction to sequences of the form
$\sum_{d\mid n}d\,g_d(n)$,
where $g_n:\mathbb{Z}_{>0}\mapsto\mathbb{Z}$ for all $n\geq1$, and determine conditions on the $g_n(\cdot)$ under which these sequences are Euler--Gauss.

\section{\texorpdfstring{$q$}{q}-Euler--Gauss Sequences}\label{sec:3}
In this section, the class of $q$-Euler--Gauss sequences is introduced. These sequences are compared with the generalisations of the $q$-Gauss sequences in the literature. To this end, conditions under which the $q$-Euler--Gauss sequences are $q$-Gauss sequences, are determined. Some new equivalents to the $q$-Gauss congruence are also obtained. The following Venn diagram represents all the containment relationships between the $q$-Euler--Gauss sequences and other generalisations of the $q$-Gauss sequences.
\begin{figure}[t]
\centering
\includegraphics[scale=0.345]{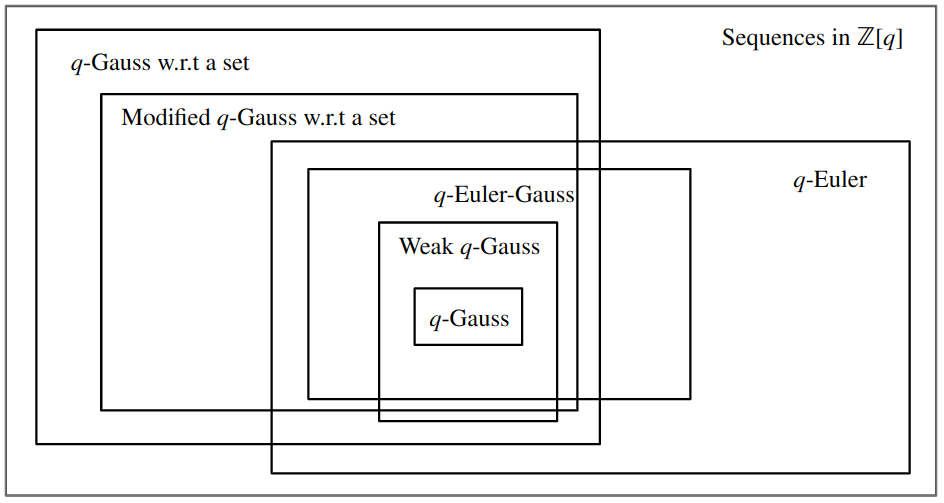}
\caption{$q$-Gauss sequences and generalisations}
\label{fig:2}
\end{figure}

\begin{definition}[$q$-Euler--Gauss Sequences]
    A polynomial sequence $(a_n(q))$ in $\mathbb{Z}[q]$ is defined to be a $q$-Euler--Gauss sequence if it satisfies 
    \begin{equation}\label{eq:17}
        \prod_{\substack{d\mid n \\ \mu(d)=1}} a_{\frac{n}{d}}(q^d) \equiv \prod_{\substack{d\mid n \\ \mu(d)=-1}} a_{\frac{n}{d}}(q^d)  \pmod{[n]_q},
    \end{equation}
    for all $n\geq1$. At $q=1$, this congruence becomes the Euler--Gauss congruence. Thus, $(a_{n}(1))$ is an Euler--Gauss sequence.
\end{definition}
Noting that at prime powers the congruence reduces to a $q$-analog of the Euler congruence, we are prompted to define the class of \textit{$q$-Euler sequences}.
\begin{definition}[$q$-Euler Sequences]
    A polynomial sequence $(a_n(q))$ in $\mathbb{Z}[q]$ is defined to be a $q$-Euler sequence if it satisfies 
    \[a_{p^r}(q)\equiv a_{p^{r-1}}(q^p)\pmod{[p^r]_q}\]
    for all primes $p$ and integers $r\geq1$.
\end{definition}
It is immediate that every $q$-Euler--Gauss sequence is a $q$-Euler sequence and that for every $q$-Euler sequence $(a_n(q))$, the integer sequence $(a_n(1))$ is an Euler sequence. We also note that by definition, not all $q$-Gauss sequences with respect to some $S$ are $q$-Euler sequences and therefore not $q$-Euler--Gauss sequences. The same holds even for integer sequences.

Before establishing that every $q$-Gauss sequence is a $q$-Euler--Gauss sequence, we prove a property of the M\"obius function that gives rise to a useful symmetry between the divisor indices of the terms appearing in the $q$-Gauss and $q$-Euler--Gauss congruences.

\begin{lemma}\label{lem:05}
    For any integers $n>1$ and $1\leq j<n$, the following multisets are equal:
    \[\{\gcd(n/d, j):d\mid n, \mu(d)=1\}=\{\gcd(n/d, j):d\mid n, \mu(d)=-1\}\]
\end{lemma}
\begin{proof}
    Let $d'$ be a divisor of $n$. It suffices to show that the multiplicities of $d'$ in both multisets are the same, or equivalently 
    \[
        \sum_{\substack{d\mid n\\\gcd(\frac{n}{d},j)=d'}}\hspace{-4mm}\mu(d)\hspace{1mm}=\hspace{1mm}0
    \]
    for every such $d'$. If $d'\nmid j$, then it does not appear in either multiset (zero multiplicity). For the divisors $d'\mid j$, let $N=n/d'$ and $J=j/d'$. Further, let $K$ denote the largest divisor $k\mid N$ such that $\operatorname{rad}(\gcd(N,J))=\operatorname{rad}(k)$ where $\operatorname{rad}(k)$ is the square-free kernel of $k$ with $\operatorname{rad}(1)=1$. Let $K=1$ if $N$ and $J$ are coprime so that $K$ is well-defined. 
    
    Since $\gcd(n/d,j)=d'$ is equivalent to the conditions $d\mid N$, $\gcd(N/d,J)=1$ and $K$ is the smallest divisor of $N$ satisfying the latter, every such $d$ may be written as $d=Ke$ for $e\mid \frac{N}{K}$. Consequently, 
    \[\sum_{\substack{d\mid N\\\gcd(\frac{N}{d},J)=1}}\hspace{-4mm}\mu(d)=\sum_{e\mid \frac{N}{K}}\mu(Ke)=\mu(K)\sum_{e\mid\frac{N}{K}}\mu(e)\]
    where second equality follows since $K$ and $N/K$ are coprime by the definition of $K$. The final sum vanishes unless $N=K$, in which case $N$ and $J$ have the same prime divisors, so either $J=N$, contradicting $j<n$ or $K$ is not square-free, so that $\mu(K)=0$. Hence the initial sum equals zero in every case.
\end{proof}

\begin{theorem}\label{thm:5}
Every $q$-Gauss sequence is a $q$-Euler--Gauss sequence.
\end{theorem}
\begin{proof}
By the definition of the $q$-Euler--Gauss congruence, it suffices to prove that for all $n>1$, $1\leq j<n$, and $n$\textsuperscript{th} primitive roots of unity $\omega_n$,
\[
\prod_{\substack{d\mid n \\ \mu(d)=1}} a_\frac{n}{d}(\omega_n^{jd})
=
\prod_{\substack{d\mid n \\ \mu(d)=-1}} a_\frac{n}{d}(\omega_n^{jd}).
\]
From the characterisation of $q$-Gauss sequences given in \eqref{eq:011}, this is equivalent to showing that
\[
\prod_{\substack{d\mid n \\ \mu(d)=1}} a_{\gcd(\frac{n}{d},j)}(1)
=
\prod_{\substack{d\mid n \\ \mu(d)=-1}} a_{\gcd(\frac{n}{d},j)}(1)
\]
for every $n>1$ and $1\leq j<n$. By the preceding lemma, the multisets of indices appearing on the left and right-hand sides are identical. Consequently, the two products consist of precisely the same factors, counted with the same multiplicities appearing in the exponent and the desired equality follows immediately.
\end{proof}

The converse of this theorem is not true. This is shown by the examples in the next subsection. The question then arises as to what restrictions are required for the converse to hold. In other words, we ask, ``Under what conditions do the $q$-Euler--Gauss sequences satisfy the $q$-Gauss congruence?" 

Before addressing this question, it may be noted that the $q$-Gauss sequences are simpler than the $q$-Euler--Gauss sequences. In fact, not only is $(a_n(1))$ a Gauss sequence if $(a_n(q))$ is a $q$-Gauss sequence, but also, for every Gauss sequence $(a_n)$, there exists a $q$-Gauss sequence $(a_n(q))$ such that $a_n(1)=a_n$ for all $n\geq1$. This follows from the characterisation \eqref{eq:022} for $q$-Gauss sequences. For the $q$-Euler--Gauss sequences, it is clear from the definition of the $q$-Euler--Gauss congruence that for every $q$-Euler--Gauss sequence $(a_n(q))$, the sequence $(a_n(1))$ is an Euler--Gauss sequence. The question then arises: ``Does there exist for every Euler--Gauss sequence  $(a_n)$, a $q$-Euler--Gauss sequence $(a_n(q))$ such that $a_n(1)=a_n$ for all $n\geq1$ ?''

Note that this question is not independent of the one asked previously. In fact, the answer to both questions is given by a single result. We show that a $q$-Euler--Gauss sequence $(a_n(q))$ is a $q$-Gauss sequence if $a_n(1)\neq0$ for all $n\geq1$. But then, it follows that if $(a_n(q))$ is not a $q$-Gauss sequence, then $a_k(1)=0$ for some $k\geq1$. Therefore, for a strictly Euler--Gauss sequence with all non-zero terms, there exists no $q$-analog of this sequence which is a $q$-Euler--Gauss sequence. To this end, we prove the following lemma.
 
\begin{lemma}\label{lem:02}
    If a sequence $(a_n(q))\subset\mathbb{Z}[q]$ satisfies the $q$-Euler--Gauss congruence at a fixed $N\geq1$ and the $q$-Gauss congruence for all $1\leq n<N$, then 
    \[
        [N]_q\text{ divides } \sum_{d\mid N} \mu(d)a_{\frac{N}{d}}(q^d)\cdot\prod_{\substack{d\mid N, d>1\\\mu(d)=1}} a_{\frac{N}{d}}(q^d)
    \]
\end{lemma}
\begin{proof}
    It suffices to show that the dividend vanishes at every root of $[N]_q$, i.e., at all the $j$\textsuperscript{th} powers of the primitive $N$\textsuperscript{th} roots of unity, $\omega_N$, for $1\leq j<N$. From the assumed $q$-Euler--Gauss congruence at $N$, the dividend is congruent modulo $[N]_q$ to 
    \[\prod_{\substack{d\mid N\\\mu(d)=-1}}a_\frac{N}{d}(q^d)+\left(\sum_{d\mid n, d>1}\mu(d)a_\frac{N}{d}(q^d)\right)\prod_{\substack{d\mid N, d>1\\\mu(d)=1}}a_\frac{N}{d}(q^d).\]
    Since $\{a_n\}_{1\leq n\leq N-1}$ satisfies the $q$-Gauss congruence property, the characterization \eqref{eq:011} implies that each term $a_\frac{N}{d}(q^d)\mid_{q=\omega_N^j}=a_{\gcd(\frac{N}{d},j)}(1)$. By Lemma \ref{lem:05}, the sum $\sum_{d\mid n, d>1}\mu(d)a_{\gcd(\frac{N}{d},j)}(1)=-a_{\gcd(N,j)}(1)$ so that the dividend becomes 
    \[\prod_{\substack{d\mid N\\\mu(d)=-1}}a_{\gcd(\frac{N}{d},j)}(1)-\prod_{\substack{d\mid N\\\mu(d)=1}}a_{\gcd(\frac{N}{d},j)}(1).\]
    As seen in the preceding theorem, this expression vanishes for all $1\leq j<N$, and the proof is complete. 
\end{proof}
 We conjecture that the divisibility condition in this lemma is a necessary condition for $q$-Euler--Gauss sequences in general, with a slight modification to the product term in the dividend. We shall elaborate on this at the end of the subsection after resolving our current discussion with the following theorem.  
\begin{theorem}\label{thm:1} 
    Let $(a_n(q))$ be a $q$-Euler--Gauss sequence. If $a_{n}(1)\neq0$ for all $n$, then it is a $q$-Gauss sequence.
\end{theorem}
\begin{proof}
We prove the contrapositive. Suppose that $(a_n(q))$ is not a $q$-Gauss sequence, and let $n\geq1$ be the least positive integer at which the $q$-Gauss congruence fails. Note that $n$ cannot be a prime power (including $1$) because $(a_n(q))$ is a $q$-Euler--Gauss sequence. Now, from Lemma \ref{lem:02}, there exist divisors $D,d'>1$ of $n$ such that $\phi_D(q)$ divides $a_\frac{n}{d'}(q^{d'})$ where $\phi_D(q)$ denotes the $D$\textsuperscript{th} cyclotomic polynomial. We show that there exists $k\geq1$ such that $a_k(1)=0$, namely, 
\[a_\frac{n}{\operatorname{lcm}(d',D)}(1)=0.\] 
As $\{a_k(q)\}_{1\leq k<n}$ satisfies the $q$-Gauss congruence at the index $n/d'<n$, by the characterisation \eqref{eq:022} of $q$-Gauss sequences there exist $\{g_k(q)\}_{1\leq k< n}\subset\mathbb{Z}[q]$ such that, 
\[
    a_\frac{n}{d'}(q)=\sum_{d\mid\frac{n}{d'}}\left[\frac{n}{d'd}\right]_{q^{d}}g_\frac{n}{d'd}(q^{d})\text{ .}
\]
Since $\phi_D(q)\mid a_\frac{n}{d'}(q^{d'})$, the value of $a_\frac{n}{d'}(q)$ at the ${d'}$\textsuperscript{th} power of the $D$\textsuperscript{th} root of unity $\omega_D=\exp{\frac{2\pi i}{D}}$ is zero. On the other hand, evaluating the above expression at $q=\omega_D^{d'}$, and writing $e=d'd$, gives
\[
a_{\frac{n}{d'}}(\omega_D^{d'})
=
\sum_{\substack{e\mid n\\d'\mid e}}
\left[\frac{n}{e}\right]_{\omega_D^e}
g_{\frac{n}{e}}(\omega_D^e).
\]
Now, only those divisors $e\mid n$ divisible by both $d'$ and $D$, equivalently by $L:=\operatorname{lcm}(d',D)$, contribute to the sum because
\[\left[\frac{n}{e}\right]_{\omega_D^e}=
\begin{cases}
    \frac{n}{e}&D\mid e,\\
    0&D\nmid e
\end{cases}
\]
Consequently, writing $e=Le'$ we obtain
\[
a_{\frac{n}{d'}}(\omega_D^{d'})
=
\sum_{e'\mid\frac{n}{L}}
\frac{n}{Le'}g_{\frac{n}{Le'}}(1),
\]
since $\omega_D^e=1$ whenever $L\mid e$. By the characterisation \eqref{eq:022}, applied at the index $n/L$, the right-hand side is precisely $a_{n/L}(1)$. Therefore
\[
0=a_{\frac{n}{d'}}(\omega_D^{d'})=a_{\frac{n}{L}}(1),
\]
as claimed.

Alternatively, $\omega_D^{d'}$ can be rewritten as $\omega_\frac{n}{d'}^\frac{n}{D}$ and we may apply the characterisation of $q$-Gauss sequences \eqref{eq:011} at the index $n/d'$ so that 
\[0=a_\frac{n}{d'}(\omega_D^{d'})=a_\frac{n}{d'}(\omega_\frac{n}{d'}^\frac{n}{D})=a_{\gcd(\frac{n}{d'},\frac{n}{D})}(1).\]
Since $\gcd(\frac{n}{d'},\frac{n}{D})$ is precisely $\frac{n}{\operatorname{lcm}(d',D)}$, we have shown that $a_\frac{n}{L}(1)=0$.
\end{proof}
Note that Theorem \ref{thm:1} highlights a limitation of the definition of the $q$-Euler--Gauss congruence. It reveals that strictly Euler--Gauss sequences with non-zero terms do not necessarily have corresponding $q$-analogs that are strictly $q$-Euler--Gauss sequences. As a consequence, it remains open to find a congruence that characterises the $q$-analogs of all the Euler--Gauss sequences.

As discussed in Section \ref{sec:1.1}, Gorodetsky provides a $q$-analog of the Pairwise Gauss congruence \eqref{eq:010}, which, unlike for the integer sequences, is not equivalent to the $q$-Gauss congruence. We provide a counterexample to demonstrate this.
\begin{example}
    Consider the sequence $(a_n(q))$ defined as, 
\[
    a_n(q)=
    \begin{cases}
        1&12\nmid n\\
        [2^k]_{q^\frac{n}{2^k3}}[3]_{q^\frac{n}{2^k3}}+a_{2^{k-1}3}\left(q^\frac{n}{2^{k-1}3}\right)&2^k3\mid n, 2^{k+1}\nmid n, k\geq2
    \end{cases}
\text{ .}\]
It satisfies the $q$-analog of the Pairwise Gauss congruence \eqref{eq:010} at all primes, but does not satisfy the $q$-Gauss congruence \eqref{eq:05} and, in fact, not even the $q$-Euler--Gauss congruence for all $n\in\{2^k3\}_{k\geq2}$. This is because, at any $n=p^rm$ for any prime factor $p$, if $12\nmid p^rm$, then $12\nmid p^{r-1}m$ so that both sides of the $q$-congruence are equal to $1$ and congruence is satisfied. If $12\mid p^rm$, then, since
\[a_{p^rm}(q)=\sum_{j=k}^2[2^j]_{q^\frac{p^rm}{2^j3}}\cdot[3]_{q^\frac{p^rm}{2^j3}}+1\equiv1\pmod{[p^r]_q}\]
and $a_{p^{r-1}m}(q^p)=1$ if $12\nmid p^{r-1}m$ or else $a_{p^{r-1}m}(q^p)$ expands like the above and is therefore congruent to $1 \pmod{[p^r]_q}$ and the $q$-congruence is satisfied. However, for $n\in\{2^k3\}_{k\geq2}$, 
\[\sum_{d\mid n}\mu(d)a_\frac{n}{d}(q^d)=[2^k]_q[3]_q\not\equiv0\pmod{[2^k3]_q}\]
and the $q$-Gauss congruence is not satisfied at such $n$.
\end{example}

The reason this $q$-congruence property does not characterise $q$-Gauss sequences is that, unlike in the case of integers where the prime powers from the prime factorisation of each $n$ appear as the divisors in the characteristic congruence, here, the $q$-integers $[p^r]_q$ need not constitute a factorisation of $[n]_q$ or a multiple of it. Seeing that all we require is a \textit{prime factorisation} of the $[n]_q$, we suggest a new congruence based on the factorisation,
\begin{equation}\label{eq:016}[n]_q=\prod_{i=1}^k[p_i^{r_i}]_{q^{\prod_{i<j\leq k}p_j^{r_j}}}
\end{equation}
where the integer $n$ has prime factorisation $\prod_{i=1}^kp_i^{r_i}$. That all the factors have no common divisors can be seen by using their cyclotomic factorisations, the formula being, 
\[[n]_q=\prod_{d\mid n,d>1}\Phi_d(q)\text{ .}\]
\begin{theorem}[Modified Pairwise $q$-Gauss Congruence]
    A sequence $(a_n(q))$ in $\mathbb{Z}[q]$ is a $q$-Gauss sequence if and only if it satisfies 
    \begin{equation}\label{eq:15}
        a_{p^rm}(q)\equiv a_{p^{r-1}m}(q^p)\pmod{[p^r]_{q^m}}
    \end{equation}
    for all primes $p$ and integers $m,r\geq1$ such that $p\nmid m$. We call this congruence as the `Modified Pairwise $q$-Gauss congruence' from here on.
\end{theorem}

\begin{proof}
    Assuming that the sequence satisfies the above congruence property, and considering an arbitrary $n$ with $n=p^rm$, $p^{r+1}\nmid n$, we see that,
    \[\sum_{d\mid n}\mu(d)a_\frac{n}{d}(q^d)=\sum_{d\mid m}\left(a_\frac{p^rm}{d}(q^d)-a_\frac{p^{r-1}m}{d}(q^{pd})\right)\mu(d)\equiv0\pmod{[p^r]_{q^m}}\text{ .}\]
 Now, this is true for every prime factor of $n$. Say $n$ has prime factorisation $n=p_1^{r_1}\cdots p_k^{r_k}$. Then since it is true that for any two coprime integers $k$ and $m$, $[k]_{q^d}$ divides $[k]_{q^m}$ if $d\mid m$, each factor of $[n]_q$ in the factorisation \eqref{eq:016} must divide the summation $\sum_{d\mid n}\mu(d)a_\frac{n}{d}(q^d)$ and the sufficient condition is proved.

 For the necessary condition, we observe that for prime powers, the $q$-Gauss congruence reduces to the congruence \eqref{eq:15}, prompting us to use this as an induction step. Then, fixing $n=p^rm$ and assuming the congruence to be true for all divisors $d\mid m$, $d<m$, we see that   
 \begin{align*}
     \sum_{d\mid n}\mu(d)a_\frac{n}{d}(q^d)&=a_{p^rm}(q)-a_{p^{r-1}m}(q^p)+\sum_{\substack{d\mid m\\d>1}}\left(a_\frac{p^rm}{d}(q^d)-a_\frac{p^{r-1}m}{d}(q^{pd})\right)\mu(d)\text{\hspace{3mm}and, }\\
     0&\equiv a_{p^rm}(q)-a_{p^{r-1}m}(q^p) \pmod{[p^r]_{q^m}}
 \end{align*}
 since $(a_n(q))$ is a $q$-Gauss sequence and $[p^r]_{q^m}\mid [n]_q$, and the induction is complete.
\end{proof}
Gorodetsky also defines the $q$-Gauss congruence with respect to $S$, \eqref{eq:017}. For integer sequences, it is known that the  Gauss congruence with respect to $S$ is equivalent to the Pairwise Gauss congruence \eqref{eq:07} satisfied for primes only in $S$. It is natural to ask if the same is true for polynomial sequences. In other words, is the $q$-Gauss congruence with respect to a set $S$ equivalent to the Modified Pairwise $q$-Gauss congruence~\eqref{eq:15}, when the latter is required only for primes in $S$? The answer is no. In fact,
\[
\text{Modified Pairwise $q$-GC for primes in $S$}
\;\implies\;
\text{$q$-GC w.r.t.\ $S$}
\;\implies\;
\text{$q$-analog of the Pairwise GC for primes in $S$}.
\]
Therefore, we propose a modified $q$-Gauss congruence with respect to a set~$S$.

\begin{theorem}[Modified $q$-Gauss Congruence with respect to a set]  
    A sequence $(a_n(q))$ in $\mathbb{Z}[q]$ satisfies the Modified Pairwise $q$-Gauss congruence for all primes from $S$ if and only if
    \begin{equation}\label{eq:032}
        \sum_{d\mid n}\mu(d)a_\frac{nm}{d}(q^d)\equiv 0\pmod{[n]_{q^m}}
    \end{equation}
    for all $n$ divisible by primes only from $S$ and all $m$ coprime to $n$.
\end{theorem}
\begin{proof}
At prime powers, this congruence reduces to the Modified Pairwise $q$-Gauss congruence and for any fixed $n$ divisible by primes only from $S$, say $n=p^rk$, $p\nmid k$ for some $p\in S$,  
\begin{align*}
    \sum_{d\mid n}\mu(d)a_\frac{nm}{d}(q^d)=\sum_{d\mid k}\left(a_\frac{p^rkm}{d}(q^d)-a_\frac{p^{r-1}km}{d}(q^{pd})\right)\mu(d)\equiv0\pmod{[p^r]_{q^{km}}}\text{.}
\end{align*}
Since this is true for all prime factors of $n$, by the factorisation \eqref{eq:016} of $[n]_q$, this congruence holds $\pmod{[n]_{q^m}}$.\\
\end{proof}

We note that when $S=\mathbb{P}$, the set of all primes, the $q$-Gauss congruence with respect to $S$,  the Modified Pairwise $q$-Gauss congruence and the Modified $q$-Gauss congruence with respect to $S$ become equivalent to the $q$-Gauss congruence property. However, the $q$-analog of the Pairwise Gauss congruence \eqref{eq:010} remains weaker than the $q$-Gauss congruence property. 

As discussed previously, Lemma \ref{lem:02} is an instance of the following generalisation. Assuming only a $q$-Euler--Gauss sequence, we conjecture that the divisibility condition here holds for all $n\geq1$, but it offers a concise proof only for $n$ with at most two distinct prime factors.
\begin{conjecture}[Divisibility criterion for $q$-Euler--Gauss sequences]\label{conj:1} 
Let $(a_n(q))$ be a $q$-Euler--Gauss sequence. Then $[n]_q$ divides $\sum_{d\mid n} \mu(d)a_{\frac{n}{d}}(q^d)\cdot A_n(q)$ for all $n\geq1$ where $A_n(q)$ is recursively defined as
\[A_n(q)=\begin{cases}
        1&n =1
        \\\prod_{\substack{d\mid n, d>1\\\mu(d)=1}}a_\frac{n}{d}(q^d)\cdot \prod_{\substack{d\mid n\\\mu(d)=-1}}A_\frac{n}{d}(q^d)&\text{otherwise}
    \end{cases}.\]
\end{conjecture}
\begin{proof}[Proof for $n$ with at most two distinct prime factors]
    We first note that by definition $A_{p^k}(q)=1$ for all primes $p$ and $k\geq1$. 
    Now, for $n$ of the form $p_1^{r_1}p_2^{r_2}$ where $p_1$, $p_2$ are primes and $r_1,r_2\geq0$, we observe that 
    \[a_{{p_1}^{r_1-1}{p_2}^{r_2-1}}(q^{p_1p_2})\prod_{1<k\leq r_1}a_{{p_1}^{r_1-k}{p_2}^{r_2-1}}(q^{{p_1}^kp_2})\prod_{1<j\leq r_1}a_{{p_1}^{r_1-1}{p_2}^{r_2-j}}(q^{p_1{p_2}^j})\text{\hspace{4mm}divides\hspace{4mm}}A_{{p_1}^{r_1}{p_2}^{r_2}}(q)\text{ .}\]
    So, we multiply this divisor term by term to $\sum_{d\mid n} \mu(d)a_{\frac{n}{d}}(q^d)$ and observe what happens $\pmod{[{p_1}^{r_1}{p_2}^{r_2}]_q}$. First, multiplying by $a_{{p_1}^{r_1-1}{p_2}^{r_2-1}}(q^{p_1p_2})$ turns the sum into \[(a_{{p_1}^{r_1}{p_2}^{r_2-1}}(q^{p_2})-a_{{p_1}^{r_1-1}{p_2}^{r_2-1}}(q^{p_1p_2}))\cdot(a_{{p_1}^{r_1-1}{p_2}^{r_2}}(q^{p_1})-a_{{p_1}^{r_1-1}{p_2}^{r_2-1}}(q^{p_1p_2}))\text{ .}\] Then, upon multiplying by $a_{{p_1}^{r_1-1}{p_2}^{r_2-2}}(q^{p_1{p_2}^2})$ and $a_{{p_1}^{r_1-2}{p_2}^{r_2-1}}(q^{{p_1}^2p_2})$, the first and second factors respectively become
    \begin{align*}
        &([p_1^{r_1}p_2^{r_2-1}]_{q^{p_2}}+a_{{p_1}^{r_1-1}{p_2}^{r_2-1}}(q^{p_1p_2})\cdot(a_{{p_1}^{r_1}{p_2}^{r_2-2}}(q^{{p_2}^2})-a_{{p_1}^{r_1-1}{p_2}^{r_2-2}}(q^{p_1{p_2}^2})))\\ &\text{ and } 
        ([p_1^{r_1-1}p_2^{r_2}]_{q^{p_1}}+a_{{p_1}^{r_1-1}{p_2}^{r_2-1}}(q^{p_1p_2})\cdot(a_{{p_1}^{r_1-2}{p_2}^{r_2}}(q^{{p_1}^2})-a_{{p_1}^{r_1-2}{p_2}^{r_2-1}}(q^{{p_1}^2p_2})))\text{ .}
    \end{align*} 
    Further multiplying by $a_{{p_1}^{r_1-1}{p_2}^{r_2-3}}(q^{p_1{p_2}^3})$ and $a_{{p_1}^{r_1-3}{p_2}^{r_2-1}}(q^{{p_1}^3p_2})$, they take the form 
    \begin{align*}
        &([p_1^{r_1}p_2^{r_2-1}]_{q^{p_2}}+f_{11}(q)\cdot([p_1^{r_1}p_2^{r_2-2}]]_{q^{{p_2}^2}}+f_{12}(q)\cdot(a_{{p_1}^{r_1}{p_2}^{r_2-3}}(q^{{p_2}^3})-a_{{p_1}^{r_1-1}{p_2}^{r_2-3}}(q^{p_1{p_2}^3}))))\\ &\text{and }([p_1^{r_1-1}p_2^{r_2}]_{q^{p_1}}+f_{21}(q)\cdot([{p_1}^{r_1-2}{p_2}^{r_2}]_{q^{{p_1}^2}}+f_{22}(q)\cdot(a_{{p_1}^{r_1-3}{p_2}^{r_2}}(q^{{p_1}^3})-a_{{p_1}^{r_1-3}{p_2}^{r_2-1}}(q^{{p_1}^3p_2}))))\text{ .}
    \end{align*}
    This continues until $a_{{p_1}^{r_1-1}}(q^{p_1{p_2}^{r_2}})$ and $a_{{p_2}^{r_2-1}}(q^{{p_1}^{r_1}p_2})$ are absorbed by these factors and they become multiples of $[{p_1}^{r_1}]_{q^{{p_2}^{r_2}}}$ and $[{p_2}^{r_2}]_{q^{{p_1}^{r_1}}}$ and their product, of $[p_1^{r_1}p_2^{r_2}]_q$.

\end{proof}
For values of $n$ with three or more distinct prime divisors, the above necessary condition for $q$-Euler--Gauss sequences has been verified computationally for $10^8$ terms.

\subsection{Examples}
In this subsection, some illustrative examples of $q$-Euler--Gauss sequences are discussed. As in the integer setting, these examples justify the strict containment relations depicted in Figure \ref{fig:2}.
\begin{example}
    Consider the sequence $(a_n(q))$ defined as,
\[a_n(q)=
    \begin{cases}
        q^n-1&12\nmid n\\
        q^\frac{n}{2^k3}&2^k3\mid n, 2^{k+1}\nmid n, k\geq2
    \end{cases}\text{.}\]    
This is a $q$-Euler--Gauss sequence but does not satisfy the $q$-Gauss congruence for $n\in\{2^k3\}_{k\geq2}$ and the $q$-analog of the Pairwise Gauss congruence at $p=2,3$. Therefore, it also does not satisfy the $q$-Gauss congruence with respect to $\{2,3\}$ and the Modified Pairwise $q$-Gauss congruence \eqref{eq:15} at $p=2,3$. This can be seen easily since at $n$ such that $12\nmid n$ or $n=2^k3$, $k\geq2$, on either side of the $q$-Euler--Gauss congruence the term $q^n-1$ appears reducing the products to $0\pmod{[n]_q}$ and for the remaining $n=2^k3^jm$, where $k\geq2$, $j\geq1$, $\gcd(m,6)=1$ and $m>1$ if $j=1$, the products become equal. Also, $a_{2^k3}(q)=q$ and $a_{2^{k-1}3}(q^2)=q^2$ breaking the $q$-analog of the Pairwise Gauss congruence at $2$ and $3$.
\end{example}

\begin{example}
    The sequence $(a_n)$ defined (recursively) as,
\[a_n(q)=
    \begin{cases}
        q-1&n=1\\
        [n]_q&n \text{ is prime}\\
        [3]_{q^{2^{k-1}}}+2&n=2^k3, k\geq2, k \text{ is even}\\
        [2]_{q^{2^{k-1}3}}+[3]_{q^{2^k}}&n=2^k3, k\geq2, k \text{ is odd}\\
        A_n(q)&\text{otherwise}
    \end{cases}\text{ .}\]
where $A_n(q)$ is the unique remainder with degree less than $n-1$ obtained by dividing the sum \[-\sum_{d\mid n, d>1}
\mu(d)a_{\frac{n}{d}}(q^d)\] by $[n]_q$.

Also a $q$-Euler--Gauss sequence, this sequence satisfies the $q$-analog of the Pairwise Gauss congruence for all primes but not the $q$-Gauss congruence for $n\in\{2^k3\}_{k\geq2}$. This can be verified by computation of the products in the $q$-Euler--Gauss congruence and the terms individually at $n=2^k3$ for all $k\geq2$, and for the remaining $n$, we know that satisfying the $q$-Gauss congruence at an $n$ is sufficient for the $q$-Euler--Gauss congruence to be satisfied at that $n$. 
Also by computation, we see that at the same $n$, the Modified Pairwise $q$-Gauss congruence fails at $2$ and $3$, the difference here being that the remainder is taken $\pmod{[2^k]_{q^3}}$ instead of $\pmod{[2^k]_{q}}$ (similarly at $p=3$). Therefore, this example further clarifies the difference between the $q$-analog of the Pairwise Gauss congruence and the Modified Pairwise $q$-Gauss congruence.
\end{example}
\begin{example}
    The sequence 
\[
    a_n(q)=
    \begin{cases}
        1&n\text{ is odd}\\
        [2]_{q^{2^{k-1}}}+1&n=2^km, k\geq1, 2\nmid m
    \end{cases}
\]
is a $q$-Euler sequence and a $q$-Gauss sequence with respect to $\{2\}$ but not a $q$-Euler--Gauss sequence. This example illustrates that there exist $q$-Euler sequences that are also $q$-Gauss sequences with respect to some $S$ that need not be $q$-Euler--Gauss sequences.
\end{example}
\begin{example}
    Another example of a $q$-Euler--Gauss sequence is the sequence
\[a_n(q)=\begin{cases}
    q+1&n=1\\
    2&n \text{ is odd}\\
    1-q^\frac{n}{2}&n=2p^k\text{ for a prime }p\\
    (2^{\Omega(n)-3}-1)q^\frac{n}{2}+(2^{\Omega(n)-3}+1)&\text{otherwise}
\end{cases}\]
where $\Omega(n)$ counts the number of prime factors of $n$, with repetitions. This is a $q$-Euler--Gauss sequence but does not satisfy the $q$-Gauss congruence at all even $n$ that are not of the form $2p^k$, where $p^k$ is a prime power. We note that from Theorem \ref{thm:2}, which shall be discussed in Section \ref{sec:5}, this sequence satisfies the $q$-Gauss congruence at all odd $n$ since $a_n(1)\neq0$ for all odd $n$. For the other examples of $q$-Euler--Gauss sequences discussed in this subsection and the next section, it may be noted that $a_1(1)=0$ and Theorem \ref{thm:2} cannot be applied.
\end{example}

With the discussion before this subsection and the above examples, Figure \ref{fig:2} is justified except for one strict containment relation: there exist $q$-Euler--Gauss sequences that are not $q$-Gauss sequences with respect to any $S$. This requires a compelling example. For this, we present a $q$-analog of the Smallest Prime Factor in the following section.

\section{Smallest Prime Factor and Greatest Prime Factor Sequences}\label{sec:4}
In this section, the class of Euler--Gauss sequences $(a_n)$ characterised by the divisor sum:
\begin{equation}\label{eq:030}
    a_n=\sum_{d\mid n}d\mu(d)g_d(n)
\end{equation}
is discussed. Here, the $g_m(\cdot)$ are integer-valued functions $g_m:\mathbb{Z}_{>0}\mapsto\mathbb{Z}$ such that, $g_1(n)=0$ for all $n\geq1$ and $g_m(n)=g_m(\text{rad}(n))$ for all $m,n\geq1$. By definition, these sequences constitute a $\mathbb{Z}$-module under pointwise addition and scalar multiplication. 
Among these sequences, there are several such as the $k$\textsuperscript{th} Smallest Prime Factor ($k$-$\operatorname{SPF}$) and $k$\textsuperscript{th} Greatest Prime Factor ($k$-$\operatorname{GPF}$) sequence, which shall be discussed further in detail. We first provide the necessary theoretical results required to show that this class of sequences indeed satisfies the Euler--Gauss congruence.
\begin{lemma}\label{lem:04}
Let $n\geq 1$ be a non-square-free integer, and let $\operatorname{rad}(n)=\prod_{p\mid n}p$
denote the radical of $n$. Then, for every divisor $d\mid n$, $\mu(d)=1$, there exists a divisor $d'\mid n$, $\mu(d')=-1$ such that
$\text{rad}\!\left(\frac{n}{d}\right)
=\text{rad}\!\left(\frac{n}{d'}\right)$.
Moreover, the divisors $d'$ may be chosen so as to define a bijection between the sets
$\{d\mid n:\mu(d)=1\}$ and
$\{d\mid n:\mu(d)=-1\}$.
\end{lemma}
\begin{proof}
Since $n$ is not square-free, at least one prime divisor of $n$ occurs with
exponent greater than one. For each divisor $d\mid n$ divisible by such a prime, let $\rho_n(d)$ denote the smallest prime divisor of $d$ that divides $n$ with exponent greater than one. Then, for each divisor $d\mid n$ with $\mu(d)=1$,
define $d'$ as follows:
\[d'=\begin{cases}
d\rho_n(n)&
p\mid\mid n\text{ for every prime divisor }p\mid d,
\text{ or }\rho_n(d)>\rho_n(n),\\[2mm]
\dfrac{d}{\rho_n(n)}&\text{otherwise}.
\end{cases}\]

By construction, $d'$ is a divisor of $n$ and satisfies $\mu(d')=-1$.
Further, since $\rho_n(n)^2\mid n$, the prime divisors of $\frac{n}{d}$ remain unchanged in $\frac{n}{d'}$, and hence their radicals agree.
To show that the correspondence between $d$ and $d'$ is one-to-one, let
$d_1$ and $d_2$ be two divisors of $n$ with $\mu(d_1)=\mu(d_2)=1$. If $d_1'$ and $d_2'$ are both obtained by the same operation (either multiplication or division by $\rho_n(n)$), then
$d_1'=d_2'$ implies $d_1=d_2$. The remaining case is $d_1\rho_n(n)=\frac{d_2}{\rho_n(n)}$
so that $d_2=d_1\rho_n(n)^2$ which contradicts the assumption that $\mu(d_2)=1$.
\end{proof}
\begin{theorem}
    Let $(a_n)$ be an integer sequence such that $a_1=0$, $p\mid a_p$ for all primes $p$ and $a_n=a_{\text{rad}(n)}$ for all $n>1$. Then $(a_n)$ is an Euler--Gauss sequence.
\end{theorem}
\begin{proof}
    For all $n\geq1$ that are not square-free, we know from the preceding lemma that for every $d\mid n$, $\mu(d)=1$, there exists a divisor $d'\mid n$ with $\mu(d')=-1$ chosen according to the bijection in the lemma such that
    \[a_\frac{n}{d}=a_{\text{rad}\left(\frac{n}{d}\right)}=a_{\text{rad}\left(\frac{n}{d'}\right)}=a_\frac{n}{d'}.\]
    Therefore, at all such $n$, $A_n^+=A_n^-$ and the Euler--Gauss congruence is satisfied. Now at all square-free $n$, if $n$ has an even number of prime factors, $\mu(n)=1$ and $a_1$ appears in $A_n^+$ so that $A_n^+=0$. Further, for every prime divisor $p$ of $n$, $\frac{n}{p}$ has an odd number of prime divisors and $a_p$ appears in $A_n^-$, but since $p\mid a_p$ for all $p$, $n\mid A_n^-$ and the Euler--Gauss congruence is satisfied at $n$. The same argument holds for square-free $n$ with an odd number of prime divisors.
\end{proof}
The sequences $(a_n)$ characterised by \eqref{eq:030} satisfy $a_1=0$ since $g_1(n)=0$ for all $n\geq1$ and $p\mid a_p$ at all primes since $a_p=-pg_p(p)$. Further, for all $n>1$, $a_{\text{rad}(n)}=\sum_{d\mid \text{rad}(n)}d\mu(d)g_d(\text{rad}(n))=a_n$, so that the following corollary is immediate.
\begin{corollary}
    The sequences $(a_n)$ characterised by the divisor sum \eqref{eq:030} are Euler--Gauss sequences. In particular, for any sequence of integer-valued functions $(f_m)$, the sequence
    \begin{equation}\label{eq:031}
    a_n=\sum_{p\mid n}pf_p(\textnormal{rad}(n)),
    \end{equation}
    where the sum is over prime divisors $p$ of $n$, is an Euler--Gauss sequence. 
\end{corollary}

The Smallest Prime Factor ($\operatorname{SPF}$) and Greatest Prime Factor ($\operatorname{GPF}$) sequences are integer sequences, (\cite[\textit{OEIS A020639}]{oeis}) and (\cite[\textit{OEIS A006530}]{oeis}) respectively, defined suitably at $1$ as
\[s_n=\begin{cases}
        0 & n=1\\
        \operatorname{SPF}(n) & \text{otherwise}
    \end{cases}
\hspace{2mm}\text{, } \hspace{10mm}g_n=\begin{cases}
        0 & n=1\\
        \operatorname{GPF}(n) & \text{otherwise}
    \end{cases}\]
where $\operatorname{SPF}(\cdot)$ denotes `smallest prime factor of' and $\operatorname{GPF}(\cdot)$, `greatest prime factor of'. One can easily see that $(s_n)$ and $(g_n)$ are defined by the prime divisor sum \eqref{eq:031} by taking $f_m(\cdot)\equiv0$ for every composite index $m$, and, for each prime index $p$, defining $f_p(\cdot)$ as follows. For the $\operatorname{SPF}$ sequence,
\[f_p(n)=
\begin{cases}
1,& p=\min\{p'\mid n:\ p'\text{ is prime}\},\\
0,& \text{otherwise}.\end{cases}\]
For the $\operatorname{GPF}$ sequence, the functions $f_p(\cdot)$ are defined similarly, with the maximum prime divisor in place of the minimum. More generally, if $f_p(n)$ is defined to take the value $1$ when $p$ is the $k$th largest prime divisor of $n$, and $0$ otherwise, then the resulting sequence is the $k$\textsuperscript{th} greatest prime factor sequence. The $k$\textsuperscript{th} smallest prime factor sequence is defined analogously.

We adopt the convention \( s_{1} = g_{1} = 0 \), since this choice is consistent with the distinct prime divisor sum \eqref{eq:031}. Moreover, these sequences satisfy Alladi’s duality identities \cite[Lemma~1]{gpfspf0}, which are expressed through the Lambert series expansions of their generating functions (see \eqref{eq:024}). This formulation again requires assuming that the first term is \(0\).

We show that the $\operatorname{SPF}$ and $\operatorname{GPF}$ sequences are counterexamples for Gauss sequences with respect to
$S$. Using the equivalents of the Gauss congruence discussed previously, we also show that these sequences satisfy Alladi's duality identities \cite[Lemma 1]{gpfspf0}. In the following section, a CSP condition for these sequences shall be presented. Both of these sequences are examples of Euler–Gauss sequences that are not Gauss sequences. In fact, neither of them satisfies the Gauss congruence for all composite squarefree numbers. Furthermore, the $\operatorname{GPF}$ sequence is a Gauss sequence only with respect to $\{2\}$, and the $\operatorname{SPF}$ sequence is not a Gauss sequence with respect to any set. Specifically, for the first $n$ terms, the $\operatorname{GPF}$ sequence does not satisfy \eqref{eq:07} for all primes $2<p\leq\frac{n}{2}$, and the $\operatorname{SPF}$ sequence for all primes $p<n$ such that there exists a prime $p'>p$ such that $pp'\leq n$. We prove these statements in the $q$-analog setting; the proofs for the integer case hold by setting $q=1$.

We note that the formal power series,  $\exp\left(\sum_{n\geq1}\frac{g_n}{n}q^n\right)$, for the $\operatorname{GPF}$ sequence belongs to $\mathbb{Z}_2[[q]]$, where $\mathbb{Z}_2$ denotes the ring of $2$-adic integers. On the other hand, we see that the $\operatorname{SPF}$ sequence belongs to $\mathbb{Q}[[q]]$. These statements are immediate from Gorodetsky's characterisation of Gauss sequences with respect to the set $\{p\}$, where $p$ denotes a prime. The Lambert series expansions of the generating functions of these sequences are, 
\[
\sum_{n\geq1}s_nq^n=\sum_{\mu(n)\neq0}(-1)^{\omega(n)+1}g_n\frac{q^n}{1-q^n}=\sum_{n\geq1}(-\mu g*\mathbf{1})_nq^n
\]
and, 
\[
\sum_{n\geq1}g_nq^n=\sum_{\mu(n)\neq0}(-1)^{\omega(n)+1}g_n\frac{q^n}{1-q^n}=\sum_{n\geq1}(-\mu s*\mathbf{1})_nq^n,
\]
where $\omega(n)$ counts the number of distinct prime divisors of $n$. In the second equality, the product $*$ denotes a Dirichlet Convolution.  
 
\begin{remark}
    Comparing the above formal power series, we have the Alladi identities: 
\begin{equation}\label{eq:024}
    s_n=(-\mu g*\mathbf{1})_n =-\sum_{d|n}\mu(d)g_{d}\text{,\hspace{4mm} and\hspace{4mm}}g_n=(-\mu s*\mathbf{1})_n=-\sum_{d|n}\mu(d)s_{d}.
\end{equation}
Note that Alladi's duality identities \cite[Lemma 1]{gpfspf0} were stated for the general case of an arithmetic function $f$ applied to the $\operatorname{SPF}$ and $\operatorname{GPF}$ functions under the assumption $f(1)=0$. In our setting, this assumption is unnecessary because we have chosen $s_{1}=g_{1}=0$. Thus, we may assert that the Alladi-type identity holds for $f(x)=x$; in particular, $f(1)=1\neq 0$.
\end{remark}

As expected, $n$ does not divide the corresponding Lambert series coefficient for all the composite square-free $n$, but does so otherwise. That the coefficients in both series expansions are `\textit{duals}' shall be verified in the discussion on the $q$-analogs of these sequences. From the Lambert series coefficients, we also obtain the generating functions:  
\[\exp\left(\sum_{m\geq1}\frac{g_m}{m}q^m\right)=\prod_{n\geq1}\frac{1}{(1-q^n)^\frac{-\mu(n)s_n}{n}}\ \in \mathbb{Z}_{2}[[q]]\] 
and, 
\[\exp\left(\sum_{m\geq1}\frac{s_m}{m}q^m\right)=\prod_{n\geq1}\frac{1}{(1-q^n)^\frac{-\mu(n)g_n}{n}}\ \in \mathbb{Q}[[q]].\]

\begin{definition}
    The $q$-analogs of the $\operatorname{SPF}$ and $\operatorname{GPF}$ sequences are defined respectively as 
 
\[S_n(q)=\begin{cases}
        0 & n=1\\
        \Phi_{\operatorname{SPF}(n)}\!\left(q^\frac{n}{\operatorname{SPF}(n)}\right) & \text{otherwise}
    \end{cases}
\hspace{2mm}\text{, }\hspace{10mm}G_n(q)=\begin{cases}
        0 & n=1\\
        \Phi_{\operatorname{GPF}(n)}\!\left(q^\frac{n}{\operatorname{GPF}(n)}\right) & \text{otherwise}\text{, }
    \end{cases}\]
 where $\Phi_k(q)$ denotes the $k$\textsuperscript{th} cyclotomic polynomial. It may be noted right away that at $q=1$, both these $q$-polynomial sequences revert to the integer $\operatorname{SPF}$ and $\operatorname{GPF}$ sequences. Also, the $q$-analogs of the $k$-$\operatorname{SPF}$ and $k$-$\operatorname{GPF}$ sequences are defined respectively as, 
 \[S_n^k(q)=\begin{cases}
        0 & n=1\\
        \Phi_{p_k(n)}(q^\frac{n}{p_k(n)}) & \text{otherwise}
    \end{cases}
\hspace{2mm}\text{, }\hspace{10mm}
G_n^k(q)=\begin{cases}
        0 & n=1\\
        \Phi_{P_k(n)}(q^\frac{n}{P_k(n)}) & \text{otherwise}\text{, }
    \end{cases}\]
where $p_k(n)$ denotes the $k$\textsuperscript{th} smallest prime factor of $n$ and $P_k(n)$ denotes the $k$\textsuperscript{th} largest prime factor of $n$. It is implicit that if $n$ has fewer than $k$ distinct prime factors, then $p_k(n)=\operatorname{GPF}(n)$ and $P_k(n)=\operatorname{SPF}(n)$.
\end{definition}

Note that the first terms $S_1(q)=0$ and $G_1(q)=0$ (for the $k$-$\operatorname{SPF}$ and $k$-$\operatorname{GPF}$ sequences as well). This is consistent with the fact that these sequences are $q$-analogs of the $\operatorname{SPF}$ and $\operatorname{GPF}$ sequences so that $S_1(1)=0=s_1$ and $G_1(1)=0=g_1$. Further, in the following section, these sequences shall be shown to generate triples that exhibit the Cyclic Sieiving Phenomenon at each index. This requires the terms $S_n(q)$ and $G_n(q)$ to have non-negative coefficients for each $n$. In particular, since $S_1(1)=G_1(1)=0$, it follows that $S_1(q)=G_1(q)=0$.
\begin{theorem}
     The sequences $(G_n(q))$ and $(S_n(q))$ are $q$-Euler--Gauss sequences and do not satisfy the $q$-Gauss congruence for all composite square-free numbers. 
\end{theorem}
\begin{proof}
       
     For non-square-free numbers written as $p^km$, $k\geq2$, $m\geq1$, where $p$ is prime, 
    \begin{equation}\label{eq:018}
     G_\frac{p^km}{d}(q)=G_\frac{p^{k-1}m}{d}(q^{pd})\text{ .}
     \end{equation}
     for all divisors $d\mid m$. From this property, one can easily see that $q$-Euler--Gauss congruence for such $n=p^km$ is satisfied. For the remaining $n$, 
     \[G_1(q^n)=0\text{ ,} \hspace{4mm}\text{ and }\hspace{4mm} \prod_{p\mid n}G_p(q^\frac{n}{p})=\prod_{p\mid n}\Phi_p(q^\frac{n}{p})\text{ .}\]
     Both of these are divisible by $[n]_q$ and always appear on opposite sides of the $q$-Euler--Gauss congruence. The same holds for $(S_n(q))$. 
     
     For composite square-free numbers, writing $n=p_1p_2\cdots p_r$, for primes $p_1<p_2<\cdots<p_r$, $r\geq1$,   
\begin{align*}
    \sum_{d\mid n}G_d\!\left(q^\frac{n}{d}\right)\mu\!\left(\frac{n}{d}\right)&=(-1)^rG_1(q^n)+\sum_{d\mid n, d>1}\Phi_{\operatorname{GPF}(d)}\!\left(q^\frac{n}{\operatorname{GPF}(d)}\right)\mu\!\left(\frac{n}{d}\right)\\
    &=0+\sum_{i=1}^r\Phi_{p_i}\!\left(q^\frac{n}{p_i}\right)\sum_{\substack{d\mid n\\\operatorname{GPF}(d)=p_i}}\mu\!\left(\frac{n}{d}\right)
\end{align*}
The coefficient of $\Phi_{p_1}(q^\frac{n}{p_1})$ is $(-1)^{r-1}$ since the divisor sum on the right-hand side runs over all the divisors of $n$ whose greatest prime factor is $p_1$ and only $p_1$ is such a divisor. For $1<i<r$, 
\[\sum_{\substack{d\mid n\\\operatorname{GPF}(d)=p_i}}\mu\!\left(\frac{n}{d}\right)=\sum_{\substack{d\mid p_1\cdots p_{i-1}}}\mu\!\left(\frac{n}{dp_i}\right)=\mu(p_{i+1}\cdots p_r)\sum_{\substack{d\mid p_1\cdots p_{i-1}}}\mu\!\left(\frac{p_1\cdots p_{i-1}}{d}\right)=0.\]
At $i=r$, the divisor sum simply becomes \[\sum_{d\mid p_1\cdots p_{r-1}}\mu\!\left(\frac{p_1\cdots p_{r-1}}{d}\right)=0.\]
Thus, the sum on the left-hand side that we started with, always equals $(-1)^{r-1}\Phi_{p_1}(q^\frac{n}{p_1})$ which is not divisible by $[n]_q$ and the Gauss congruence fails. 

The same proof works for $(S_n)$ with the only exception that the coefficient of $\Phi_{p_r}(q^\frac{n}{p_r})$ doesn't vanish and is equal to $(-1)^{r-1}$. 
\end{proof}
We note that by setting $q=1$, the theorem immediately implies that the $\operatorname{SPF}$ and $\operatorname{GPF}$ sequences, $(s_n)$ and $(g_n)$, are Euler--Gauss sequences. It may also be observed that at $q=1$, the computation in the preceding proof reduces to 
\[\sum_{d\mid p_1p_2\cdots p_r}g_\frac{p_1p_2\cdots p_r}{d}\mu(d)=(-1)^{r-1}p_1\]
which proves the relationship between the Lambert series coefficients at composite square-free $n$ for the sequence $(g_n)$. For the remaining $n$, from \eqref{eq:018} it follows that the coefficient evaluates to $0$. The same holds for the sequence $(s_n)$.

\begin{theorem}
    The sequence $(G_n(q))$ satisfies the $q$-Gauss congruence only with respect to the set $\{2\}$ and the sequence $(S_n(q))$ does not satisfy the $q$-Gauss congruence with respect to any set $S$. 
\end{theorem}
\begin{proof}
     Let $p$ be a prime. Consider the $q$-analog of the Pairwise Gauss congruence \eqref{eq:010} at $p$ for both these sequences.
     
     For the sequence $(G_n(q))$, if $p>2$, then $G_{2p}(q)=\Phi_p(q^2)$ and $G_2(q^p)=\Phi_2(q^p)$ so that $[p]_q$ divides the former but not the latter. Hence the $q$-analog of the Pairwise Gauss congruence fails at $p$. However, at $p=2$, for all $k,m\geq1$, $G_{2^km}(q)=G_{2^{k-1}m}(q^2)$ by the definition of $(G_n(q))$ so that the $q$-Gauss congruence is satisfied with respect to $\{2\}$.
     
     For the sequence $(S_n(q))$, there exists a prime $p'>p$ so that $S_{pp'}(q)=\Phi_p(q^{p'})$ and $S_{p'}(q^p)=\Phi_{p'}(q^p)$. Since $[p]_q$ divides the former but not the latter, the $q$-analog of the Pairwise Gauss congruence fails at $p$. As $p$ is an arbitrary prime, the $q$-Gauss congruence is not satisfied with respect to any $S$. 
\end{proof}
Specifically, for the first $N$ terms of the sequence $(G_n(q))$, the $q$-analog of the Pairwise Gauss congruence fails for every prime $2<p\leq\frac{N}{2}$. Indeed, for such a prime $p$, we have $2p\leq N$, and the proof above shows that the $q$-analog of the Pairwise Gauss congruence fails at $p$.

Similarly, for the first $N$ terms of $(S_n(q))$, the $q$-analog of the Pairwise Gauss congruence fails for every prime $p<N$ for which there exists a prime $p'>p$ satisfying $pp'\leq N$.

The above discussion yields two conclusions: not all Euler--Gauss sequences are Gauss sequences with respect to some $S$, and not all $q$-Euler--Gauss sequences are $q$-Gauss sequences with respect to some $S$.

\section{Cyclic Sieving Phenomenon}\label{sec:5}
In Gorodetsky's \cite{goro_q-gauss} and Gossow's \cite{fern_gauss} works, the characterising condition for Gauss sequences,  
\begin{equation}\label{eq:025}
    a_n(\omega_n^j)=a_{\gcd(n,j)}(1),
\end{equation}
discussed in Section \ref{sec:1.1} links the $q$-Gauss sequences to the Cyclic Sieving Phenomenon (CSP). For the $q$-Euler--Gauss sequences, one can see that this condition holds at prime powers. Further, from Theorem \ref{thm:1}, it is established that the $q$-Euler--Gauss sequences are $q$-Gauss sequences provided that all the terms are non-zero at $q=1$. The next result describes the conditions under which $a_{N}(q)$, for a fixed $N \geq 1$, satisfies the $q$-Gauss congruence, assuming that the sequence $(a_{n}(q))$ is $q$-Euler--Gauss. This observation allows us to apply the CSP condition for the $q$-Gauss sequences,
\begin{equation}\label{eq:026}
    |X_n^j| = |X_{\gcd(n,j)}|,
\end{equation}
selectively to a different sequence of sets, enabling them to form CSP triples using a $q$-Euler--Gauss sequence.

\begin{theorem}\label{thm:2}   
    Let $(a_n(q))$ be a $q$-Euler--Gauss sequence. For an integer $N>1$, if $a_d(1)\neq0$ for all proper divisors $d\mid N$ such that $\frac{N}{d}$ is not a prime power, then $(a_n(q))$ satisfies the $q$-Gauss congruence at $N$.
\end{theorem}
\begin{proof}
For any integer $n>1$, let us define
\[E_n=\{d: d\mid n, d<n, \frac{n}{d}\text{ is not a prime power}\}.\]
If $n$ is a prime power, then $E_n$ is empty. In this case, the statement of the theorem vacuously holds, which is consistent with the fact that $q$-Gauss congruence and $q$-Euler--Gauss congruence coincide at prime powers.

We first establish two elementary properties of the sets $E_n$.
\begin{enumerate}[label=(\roman*.)]
    \item For every proper divisor $d\mid n$, \[E_\frac{n}{d}\subseteq E_n.\]  
    Indeed, the claim is immediate when $n$ is a prime power. Otherwise, let $d'\in E_{n/d}$. Then $\frac{n}{dd'}$ is not a prime power. Since $\frac{n}{d'}=d\cdot\frac{n}{dd'}$, it follows that $\frac{n}{d'}$ is also not a prime power, and hence $d'\in E_n$.
    \item \[E_n=\hspace{-2.5mm}\bigcup\limits_{\substack{d\mid n, d>1\\\mu(d)=1}}\hspace{-2mm}D_\frac{n}{d}\] 
    where $D_k$ denotes the set of divisors of $k$ (including $1$ and $k$).
    
    Again, the claim is immediate when $n$ is a prime power. Suppose then that $n$ is not a prime power. If $d\in E_n$, then $\frac{n}{d}$ has at least two distinct prime divisors. Consequently, there exists a square-free divisor $d'>1$ of $\frac{n}{d}$ with $\mu(d')=1$. Hence $d\mid \frac{n}{d'}$, so that $d\in D_{n/d'}$. The reverse inclusion follows by reversing the argument.
\end{enumerate}

\vspace{-2mm}
We now proceed to prove the theorem by applying induction on the divisors of $N$. The base case is that of prime $N$ and has already been addressed. Let us assume then that for all proper divisors $d\mid N$, the $q$-Gauss congruence holds at $d$ whenever $a_{d'}(1)\neq0$ for all $d'\in E_d$. Now since the $q$-Euler--Gauss congruence is satisfied at $N$, 
\[\prod_{\substack{d\mid N \\ \mu(d)=1}} a_\frac{N}{d}(\omega_N^{jd})
=\prod_{\substack{d\mid N \\ \mu(d)=-1}} a_\frac{N}{d}(\omega_N^{jd})
\]
where $\omega_N$ is any $N$\textsuperscript{th} primitive root of unity and $1\leq j<N$. By Property (i.), the assumption $a_{d'}(1)\neq 0$ for all $d'\in E_N$ implies that $a_{d''}(1)\neq 0$ for all $d''\in E_\frac{N}{d}$ for every proper divisor $d\mid N$. In particular, this holds for every square-free divisor $d\mid N$, which are precisely the indices occurring in the above products. Hence, by the induction hypothesis the $q$-Gauss congruence holds for each $\frac{N}{d}$, and therefore
$a_\frac{N}{d}(\omega_N^{jd})=a_\frac{N}{d}(\omega_\frac{N}{d}^{\,j})=a_{\gcd(\frac{N}{d},j)}(1)$,
where the last equality follows from the characterisation \eqref{eq:011} of $q$-Gauss sequences. Substituting this into the $q$-Euler--Gauss congruence gives
\[
a_N(\omega_N^j)\prod_{\substack{d\mid N, d>1\\\mu(d)=1}}a_{\gcd(\frac{N}{d},j)}(1)=\prod_{\substack{d\mid N\\\mu(d)=-1}}a_{\gcd(\frac{N}{d},j)}(1)
\]
which by Lemma \ref{lem:05} equals
\[
\prod_{\substack{d\mid N\\\mu(d)=1}}a_{\gcd(\frac{N}{d},j)}(1).
\]
Finally, every index $\gcd\left(\frac{N}{d},j\right)$ appearing in the above product belongs to $D_{\frac{N}{d}}$ for some divisor $d>1$ with $\mu(d)=1$. Property (ii.) therefore implies that each such index belongs to $E_N$. 

By hypothesis, all corresponding terms are non-zero and cancel, leaving $a_N(\omega_N^j)=a_{\gcd(N,j)}(1)$. The characterisation \eqref{eq:011} now implies that the $q$-Gauss congruence holds at $N$ and completes the induction. 
\end{proof}
Based on this result, one can collect together the terms in a $q$-Euler--Gauss sequence that satisfy the $q$-Gauss congruence, and generate new examples of triples that exhibit the CSP. This leads us to the CSP for $q$-Euler--Gauss sequences as a corollary.

\begin{corollary}\label{coro:1} 
Let $(X_n)$ be a sequence of sets and $(a_n(q))$ be a $q$-Euler--Gauss sequence such that each $a_n(q)$ has non-negative integer coefficients. Then, for all $n>1$ such that $a_d(1)\neq0$ for all proper divisors $d\mid n$ such that $\frac{n}{d}$ is not a prime power, the triple $(X_n,\mathbb{Z}_n,a_n(q))$ satisfies the Cyclic Sieving Phenomenon if $|X_d|=a_d(1)$ for all $d\mid n$ and $\mathbb{Z}_n$ acts on $X_n$ in such a way that, 
\[|X_n^j|=|X_{\gcd(n,j)}| \text{ for all }j\in\mathbb{Z}_n\]
where $X_n^j$ is the fixed point set of $j$.
\end{corollary}

The $q$-analogs of the $\operatorname{SPF}$ and $\operatorname{GPF}$ sequences, being $q$-Euler--Gauss sequences, naturally satisfy the CSP as discussed above. However, for these sequences, another CSP condition arises for the sequence of sets $(X_n)$ that may easily be contrasted with the standard CSP condition \eqref{eq:026} for the $q$-Gauss sequences. This is due to the following behaviour of these sequences when evaluated at roots of unity. 
\begin{theorem}\label{thm:3}
    The sequences $(S_n(q))$ and $(G_n(q))$ when evaluated at the $n$\textsuperscript{th} roots of unity given by $\omega_n^j$, $j\geq1$ where $\omega_n$ is a primitive root, satisfy,
    \[S_n(\omega_n^j)=S_{\gcd(s_n,j)}(1)\text{ ,\hspace{4mm}and\hspace{5mm}} G_n(\omega_n^j)=G_{\gcd(g_n,j)}(1)\]
for all integers $n>1$. 
\end{theorem}
\begin{proof}
    From the definitions of the sequences, they satisfy
    \[S_n(\omega_n^j)=\begin{cases}
        s_n&\text{for }s_n\mid j\\
        0&\text{otherwise}
\end{cases}\text{ ,\hspace{4mm}and\hspace{5mm}} G_n(\omega_n^j)=\begin{cases}
        g_n&\text{for }g_n\mid j\\
        0&\text{otherwise}
    \end{cases}\text{ .}\]
    For $S_n(q)$, we see that $S_n(\omega_n^j)=\Phi_{s_n}(\omega_n^\frac{nj}{s_n})$, which is equal to zero unless $s_n\mid j.$ When $s_n\mid j$, $S_n(\omega_n^j)=s_n$. The same holds for the sequence $(G_n(q))$.
\end{proof}
From this, it remains only to find a sequence of sets enumerated by these sequences, and the following corollary is immediate. 

\begin{corollary}
    The $n$\textsuperscript{th} term of the sequence $(S_n(q))$, the set $X_n$ from sequence of sets $(X_n)$ such that $X_1=\emptyset$ and $|X_p|=p$ for all primes $p$, and the addition modulo $n$ group,  $\mathbb{Z}_n$ acting on $X_n$ in such a way that, 
    \[|X_n^j|=|X_{\gcd(s_n,j)}|\] for all $j\in\mathbb{Z}_n$, constitute a CSP triple $(X_n,\mathbb{Z}_n,S_n(q))$ for all $n\geq1$. The same result holds for $(G_n(q))$ with $s_n$ replaced by $g_n$ in the above equality.
\end{corollary}
Note that our choice of $S_1(q)=G_1(q)=0$ is consistent with the CSP exhibited by the triples $(X_n,\mathbb{Z}_n,S_n(q))$ and $(X_n,\mathbb{Z}_n,G_n(q))$ at $n=1$. In the preceding result, the set $X_1$ is required to be empty which forces the polynomials $S_1(q)$ and $G_1(q)$ to be identically zero. 

\section*{Acknowledgement} 
We thank the anonymous referee for carefully reading the manuscript and providing detailed and insightful comments. We especially acknowledge the suggestion on inclusion of the coprimality condition (i.) in Definition \ref{defn:strongEG}, and the correction of Example \ref{ex:4}. The second author is supported by the SERB-MATRICS project (MTR/2023/000705) from the Department of Science and Technology, India, for this work.

\bibliographystyle{plainurl}
\bibliography{references}

@article{armin_constant,
title = {On the representability of sequences as constant terms},
journal = {Journal of Number Theory},
volume = {253},
pages = {235-256},
year = {2023},
issn = {0022-314X},
doi = {https://doi.org/10.1016/j.jnt.2023.06.015},
url = {https://www.sciencedirect.com/science/article/pii/S0022314X23001373},
author = {Alin Bostan and Armin Straub and Sergey Yurkevich}
}

@article{goro_q-gauss,
    author = {Gorodetsky, Ofir},
    title = {q-Congruences, with applications to supercongruences and the cyclic sieving phenomenon},
    journal = {International Journal of Number Theory},
    volume = {15},
    number = {09},
    pages = {1919-1968},
    year = {2019},
    doi = {10.1142/S1793042119501069},
}

@article{dold,
    author = "Graff, Grzegorz and Gulgowski, Jacek and Lebiedź, Małgorzata",
    year = "2022",
    title = "Generalized Dold Sequences on Partially-Ordered Sets",
    volume = "29",
    journal = "The Electronic Journal of Combinatorics",
    doi = "10.37236/10544"
}

@misc{fern_gauss,
      title={Polynomial and combinatorial analogues of Gauss congruence}, 
      author={Fern Gossow},
      year={2024},
      eprint={2410.05678},
      archivePrefix={arXiv},
      primaryClass={math.CO},
      url={https://arxiv.org/abs/2410.05678}, 
}

@article{minton_gauss,
author = {Minton, Gregory},
year = {2014},
pages = {},
title = {Linear recurrence sequences satisfying congruence conditions},
volume = {142},
journal = {Proceedings of the American Mathematical Society},
doi = {10.1090/S0002-9939-2014-12168-X}
}

@book{stanley_2,
  author = {Stanley, Richard P.},
  title = {Enumerative Combinatorics, Volume 2},
  year = {1999},
  publisher = {Cambridge University Press},
  isbn = {978-0-521-56069-6},
  edition = {2nd},
}

@article{csp,
title = {The cyclic sieving phenomenon},
journal = {Journal of Combinatorial Theory, Series A},
volume = {108},
number = {1},
pages = {17-50},
year = {2004},
issn = {0097-3165},
doi = {https://doi.org/10.1016/j.jcta.2004.04.009},
url = {https://www.sciencedirect.com/science/article/pii/S0097316504000822},
author = {V. Reiner and D. Stanton and D. White},
}

@article{Smyth,
author = {C. J. Smyth},
title = {A Coloring Proof of a Generalisation of Fermat's Little Theorem},
journal = {The American Mathematical Monthly},
volume = {93},
number = {6},
pages = {469--471},
year = {1986},
publisher = {Taylor \& Francis},
doi = {10.1080/00029890.1986.11971858},
URL ={https://doi.org/10.1080/00029890.1986.11971858}}

@article{zarelua,
author = {Zarelua, A.},
year = {2008},
pages = {78-98},
title = {On congruences for the traces of powers of some matrices},
volume = {263},
journal = {Proceedings of the Steklov Institute of Mathematics},
doi = {10.1134/S008154380804007X}
}

@book{dickson,
  title={History of the Theory of Numbers},
  author={Dickson, Leonard Eugene},
  volume={1},
  year={1919},
  publisher={Carnegie Institution of Washington},
  address={Washington},
  series={History of the Theory of Numbers},
  url={https://doi.org/10.5962/t.174912}
}

@article{dold_survey,
   title={Dold sequences, periodic points, and dynamics},
   volume={53},
   ISSN={1469-2120},
   url={http://dx.doi.org/10.1112/blms.12531},
   DOI={10.1112/blms.12531},
   number={5},
   journal={Bulletin of the London Mathematical Society},
   publisher={Wiley},
   author={Byszewski, Jakub and Graff, Grzegorz and Ward, Thomas},
   year={2021},
   pages={1263–1298} }

@article{generalized_fermat,
title = {Generalized Fermat, double Fermat and Newton sequences},
journal = {Journal of Number Theory},
volume = {98},
number = {1},
pages = {172-183},
year = {2003},
issn = {0022-314X},
doi = {https://doi.org/10.1016/S0022-314X(02)00025-2},
url = {https://www.sciencedirect.com/science/article/pii/S0022314X02000252},
author = {Bau-Sen Du and Sen-Shan Huang and Ming-Chia Li}
}

@article{csp-lyndon,
   title={The Cyclic Sieving Phenomenon on Circular Dyck Paths},
   volume={26},
   ISSN={1077-8926},
   url={http://dx.doi.org/10.37236/8720},
   DOI={10.37236/8720},
   number={4},
   journal={The Electronic Journal of Combinatorics},
   publisher={The Electronic Journal of Combinatorics},
   author={Alexandersson, Per and Linusson, Svante and Potka, Samu},
   year={2019},
}

@article{gpfspf0,
title = {Duality between prime factors and an application to the prime number theorem for arithmetic progressions},
journal = {Journal of Number Theory},
volume = {9},
number = {4},
pages = {436-451},
year = {1977},
issn = {0022-314X},
doi = {https://doi.org/10.1016/0022-314X(77)90005-1},
url = {https://www.sciencedirect.com/science/article/pii/0022314X77900051},
author = {Krishnaswami Alladi}
}

@misc{gpfspf1,
author = {Alladi, Krishnaswami and Johnson, Jason},
year = {2024},
pages = {},
title = {Duality between prime factors and the Prime Number Theorem for Arithmetic Progressions -- II},
doi = {10.48550/arXiv.2410.18259},
url = {https://arxiv.org/abs/2410.18259},
}

@book{antapostol,
  author    = {Apostol, Tom M.},
  title     = {Introduction to Analytic Number Theory},
  series    = {Undergraduate Texts in Mathematics},
  publisher = {Springer},
  address   = {New York, NY},
  year      = {1976},
  doi       = {10.1007/978-1-4757-5579-4},
  isbn      = {978-0-387-90163-3},
  edition   = {1},
  pages     = {xii+340},
  url       = {https://doi.org/10.1007/978-1-4757-5579-4}
}

@article{potencyofn,
    author = {MacMahon, P. A.},
    title = {Dirichlet Series and the Theory of Partitions},
    journal = {Proceedings of the London Mathematical Society},
    volume = {s2-22},
    number = {1},
    pages = {404-411},
    year = {1924},
    month = {01},
    issn = {0024-6115},
    doi = {10.1112/plms/s2-22.1.404},
    url = {https://doi.org/10.1112/plms/s2-22.1.404},
    eprint = {https://academic.oup.com/plms/article-pdf/s2-22/1/404/4372558/s2-22-1-404.pdf},
}

@article{alladiintegerlog,
author = {K. Alladi and P. Erdős},
title = {{On an additive arithmetic function.}},
volume = {71},
journal = {Pacific Journal of Mathematics},
number = {2},
publisher = {Pacific Journal of Mathematics, A Non-profit Corporation},
pages = {275 -- 294},
year = {1977},
}

@misc{oeis,
  author = {{OEIS Foundation Inc.}},
  Title = {The {O}n-{L}ine {E}ncyclopedia of {I}nteger {S}equences},
  url    = {https://oeis.org},
}

@article{Wojcik_mobius_replace_gauss,
  author    = {Klaudiusz W{\'o}jcik},
  title     = {Newton Sequences and Dirichlet Convolution},
  journal   = {Integers},
  volume    = {21},
  pages     = {Paper No. A63},
  year      = {2021},
  publisher = {egang.org},
  url       = {http://math.colgate.edu/~integers/v63/v63.pdf}
}

@misc{dwork_congruence,
      title={Dwork Congruences and reflexive Polytopes}, 
      author={Kira Samol and Duco van Straten},
      year={2009},
      eprint={0911.0797},
      archivePrefix={arXiv},
      primaryClass={math.AG},
      url={https://arxiv.org/abs/0911.0797}, 
}
\end{document}